\documentclass[12pt]{article}
\newcommand\SHORTTITLE{Adhesive contact at mixed mode \& numerics}
\def\version{30.7.2013}\def\users{world}
\usepackage{amsmath,amssymb,amsthm,mathrsfs,epsf,a4,color,graphicx,eucal}
\usepackage{cite}  

\usepackage{psfrag}
\newcounter{myfigure}
\newenvironment{my-picture}[3]{\refstepcounter{myfigure}\label{#3}\setlength{\unitlength}{\textwidth}\begin{picture}(#1,#2)}{\end{picture}}

\usepackage{ifthen}\usepackage{fancyhdr}
\ifthenelse{\equal{\users}{world}}{}
{\pagestyle{fancy}\headheight=28pt
\definecolor{grey}{rgb}{0.6,0.6,0.6}
\rhead{\color{grey}\SHORTTITLE\\by M.Kruzik \& C.G.Panagiotopoulos \& T.Roubicek}
\chead{} \lhead{\color{grey}Version \version, file:
\jobname.tex\\compiled:
\number\day.\number\month.\number\year\ at
\the\hour:\ifnum\minute<10 0\fi\the\minute\ h }

\newcount\hour \newcount\minute
\hour=\time \divide \hour by 60 \minute=\time \loop \ifnum \minute >
59 \advance \minute by -60 \repeat
}
\usepackage[normalem]{ulem}
\definecolor{brown}{rgb}{0.5,0,0}
\ifthenelse{\equal{\users}{world}}{

    \newcommand{\COMMENT}[1]{}
    \newcommand{\DELETE}[1]{}
    
    }
{
 
 \newcommand{\DELETE}[1]{{\color{brown}\sout{#1}\color{black}}}
 \newcommand{\COMMENT}[1]{{\color{red}\uuline{#1}\color{black}}}

\usepackage[notcite,notref]{showkeys}
}
\definecolor{ddmagenta}{rgb}{0.7,0,0.9}
\definecolor{ddcyan}{rgb}{0,0.2,1.0}
\definecolor{dred}{rgb}{.8,0,0}

\newtheorem{theorem}{Theorem}[section]

\newtheorem{definition}[theorem]{Definition}

\newtheorem{proposition}[theorem]{Proposition}

\newtheorem{remark}[theorem]{Remark}

\numberwithin{equation}{section}

\newcommand{\ITEM}[2]{\parbox[t]{.05\textwidth}{#1}\hfill\parbox[t]{.95\textwidth}{#2}\vspace*{.8mm}}

\newcommand{\DDD}[3]{\begin{array}[t]{c}#1\vspace*{-1em}\\_{#2}\vspace*{-.5em}\\_{#3}\end{array}}
\newcommand{\ddd}[3]{\DDD{\begin{array}[t]{c}\underbrace{#1}\vspace*{.6em}\end{array}}{\text{\footnotesize #2}}{\text{\footnotesize #3}}}

\newcommand\R{\mathbb R}
\newcommand\N{\mathbb N}
\newcommand\Q{\mathcal Q}

\newcommand\DT[1]{\mathchoice
                 {{\buildrel{\hspace*{.1em}\text{\LARGE.}}\over{#1}}}
                 {{\buildrel{\hspace*{.1em}\text{\Large.}}\over{#1}}}
                 {{\buildrel{\hspace*{.1em}\text{\large.}}\over{#1}}}
                 {{\buildrel{\hspace*{.1em}\text{\large.}}\over{#1}}}}
\newcommand\DDT[1]{\mathchoice
   {{\buildrel{\hspace*{.1em}\text{\LARGE.\hspace*{-.1em}.}}\over{#1}}}
   {{\buildrel{\hspace*{.1em}\text{\Large.\hspace*{-.1em}.}}\over{#1}}}
   {{\buildrel{\hspace*{.1em}\text{\large.\hspace*{-.1em}.}}\over{#1}}}
   {{\buildrel{\hspace*{.1em}\text{\large.\hspace*{-.1em}.}}\over{#1}}}}

\newcommand\JUMP[1]{\mathchoice
                   {\big[\hspace*{-.3em}\big[#1\big]\hspace*{-.3em}\big]}
                   {[\hspace*{-.15em}[#1]\hspace*{-.15em}]}
                   {[\![#1]\!]}
                   {[\![#1]\!]}}
\newcommand\bbA{\mathbb A}
\newcommand\bbC{\mathbb C}
\newcommand\bbD{\mathbb D}

\newcommand\bbB{\mathbb B}
\newcommand\bbH{\mathbb H}

\renewcommand\d{\mathrm d}

\newcommand\overlineGC{\overline{\Gamma}_{\mbox{\tiny\rm C}}}
\newcommand\overlineSC{\overline{\Sigma}_{\mbox{\tiny\rm C}}}

\newcommand{\FRM}{F}

\newcommand{\fRM}{f}

\newcommand{\calD}{\mathcal{R}}
\newcommand{\calE}{\mathcal{E}}

\newcommand{\aein}{\text{a.e. in}}

\newcommand{\testu}{v}

\newcommand{\GD}{\mathchoice
                  {\Gamma_{\hspace*{-.15em}\mbox{\tiny\rm D}}}
                  {\Gamma_{\hspace*{-.15em}\mbox{\tiny\rm D}}}
                  {\Gamma_{\hspace*{-.1em}\mbox{\tiny\rm D}}}
                  {\Gamma_{\hspace*{-.05em}\mbox{\tiny\rm D}}}}
\newcommand{\GN}{\Gamma_{\hspace*{-.15em}\mbox{\tiny\rm N}}}
\newcommand\GC{\mathchoice{\Gamma_{\hspace*{-.15em}\mbox{\tiny\rm C}}}
                          {\Gamma_{\hspace*{-.15em}\mbox{\tiny\rm C}}}
                          {\Gamma_{\hspace*{-.05em}\mbox{\tiny\rm C}}}
                          {\Gamma_{\hspace*{-.05em}\mbox{\tiny\rm C}}}}
\newcommand{\Sdir}{\Sigma_{\mbox{\tiny\rm D}}}
\newcommand{\Snew}{\Sigma_{\mbox{\tiny\rm N}}}
\newcommand\SC{\Sigma_{\mbox{\tiny\rm C}}}
\newcommand{\psiG}{\psi_{\scriptscriptstyle\rm G}}
\newcommand{\dela}{\bbA}

\newcommand{\ind}{I}
\newcommand{\dd}{\mathrm{d}}

\newcommand{\pwc}[2]{\overline{#1}_{\kern-1pt#2}}
\newcommand{\upwc}[2]{\underline{#1}_{\kern-1pt#2}}

\newcommand{\pwl}[2]{#1_{\kern-1pt#2}}

\newcommand{\weaksto}{\stackrel{*}{\rightharpoonup}}

\newcommand{\uk}{u_{\tau}^k}

\newcommand{\zk}{z_{\tau}^k}

\newcommand{\ukm}{u_{\tau}^{k-1}}

\newcommand{\Omegaone}{\Omega_+}
\newcommand{\Omegatwo}{\Omega_-}
\newcommand{\barOmegaone}{\bar\Omega_+}
\newcommand{\barOmegatwo}{\bar\Omega_-}

\newcommand{\norm}{\vec{n}_{\mbox{\tiny\rm C}}^{}}

\newcommand{\calR}{\mathcal{R}}

\begin{document}

\noindent{\LARGE\bf
Quasistatic adhesive contact delaminating
\\[.2em]
in mixed mode and its numerical treatment}

\bigskip\bigskip

\noindent{\large\sc Martin Kru\v z\'\i k}\footnote{corresponding author, email: kruzik@utia.cas.cz}\\
{\it Institute of Information Theory and Automation of the ASCR, Pod vod\'{a}renskou v\v{e}\v{z}\'{\i} 4, CZ--182~08 Praha~8, Czech Republic}
and
{\it Faculty of Civil Engineering, Czech Technical University, Th\'{a}kurova 7, CZ--166~29 Praha~6, Czech Republic}

\medskip

\noindent{\large\sc Christos G.\ Panagiotopoulos}\\
{\it Group of Elasticity and Strength of Materials, 
Dept.~of Continuum 
Mech., School of Engineering,
University of Seville, Camino de los Descubrimientos s/n, E-41092 Sevilla, Spain}

\medskip

\noindent{\large\sc Tom\'a\v s Roub\'\i\v cek}\\
{\it Mathematical Institute, Charles University, Sokolovsk\'a 83,
CZ--186~75~Praha~8,  
} and
{\it Institute of Thermomechanics of the ASCR, Dolej\v skova~5,
CZ--182~08 Praha 8, Czech Republic}



\bigskip



\begin{center}\begin{minipage}[t]{16.5cm}

\baselineskip=12pt

{\small \noindent {\it Abstract.} 
An adhesive unilateral contact between visco-elastic bodies at
small strains and in a Kelvin-Voigt rheology is scrutinized,
neglecting inertia. The flow-rule
for debonding the adhesive is considered rate independent,
unidirectional, and non-associative due to dependence on the mixity
of modes of delamination, namely Mode I (opening) needs
(=dissipates) less energy than Mode II (shearing). 
Such mode-mixity dependence of delamination is a very
pronounced (and experimentally confirmed) phenomenon typically
considered in engineering models.

An efficient semi-implicit-in-time FEM discretization leading 
to recursive quadratic mathematical programs is devised.
Its convergence and thus the existence of weak solutions is proved.
Computational experiments implemented by BEM illustrate the modeling 
aspects and the numerical efficiency of the discretization.

\medskip


\medskip

}

\end{minipage}
\end{center}

\vspace{1cm}

\section{Introduction}\label{sec-intro}


In this article, we focus on adhesive contacts, which is 
a part of {\it nonlinear contact mechanics}
with  numerous
practical applications. In particular, we concentrate on the modeling, analysis, 
and computations of an inelastic process called delamination 
(or debonding) of elastic bodies glued together along a
prescribed delamination interfaces. 
On a microscopic level, some
macromolecular chains of the adhesive may break upon loading and we assume
that they can never be glued back, i.e., no ``healing'' is possible.
This makes the process {\it unidirectional} and  irreversible. On
the glued interface, we consider the {\it delamination process} as
{\it rate-independent} and, in the bulk, we also take into account
rate-dependent {\it viscous effects}. 
An important feature  appearing in  engineering modeling
(and so far mostly omitted in the mathematical literature), is  the
dependence of this process  on the  modes under which it proceeds.
Indeed, Mode I (=opening) usually dissipates much less energy than
Mode II (=shearing). The difference may be up to  hundreds of
percents, cf.\
\cite{BanAsk00NFCI,LieCha92ASIF,Mant08DRLM,SwLiLo99ITAM}. Moreover,
the delamination process seldom occurs in such pure modes
and, in reality, the {\it mixed mode} appears more frequently. The substantial
difference in the dissipation in various modes is explained either by some 
roughness of the glued interface (to be overcome in Mode II but not in 
Mode I, cf.\ \cite{ERDC90FEBI}) or by some plastification caused by shear 
in Mode II (but not by mere tension in Mode I) 
before the delamination itself happens,
cf.~\cite{LieCha92ASIF,TveHut93IPMM}.

This article  addresses  a standard engineering model 
which was rigorously analyzed already in \cite{rr+tr2} even in
a full thermodynamical context but exploiting the concept of 
non-simple materials (see e.g.~\cite{toupin}) which would be much more 
demanding to be implemented computationally. Although computational 
simulations are routinely performed in engineering and successfully used 
in applications even 
in simple materials, where the strain energy depends only of the first 
gradient of the deformation, cf.\ e.g.\ 
\cite{BanAsk00NFCI,CaDaMo03NSMM,HutSuo92MMCL,TMGCP10ACTA,TMGP11BEMA} 
and references therein, the rigorous convergence/existence analysis 
is not at disposal. Such computations and the models themselves are thus
completely unjustified. 
Here we concentrate on an {\it isothermal} situation in {\it simple visco-elastic 
materials}, and emphasize numerical aspects including the {\it efficient 
computational feasibility} of the model with a {\it guaranteed
numerical stability and convergence}. 


Let us just highlight  main ingredients of the
model formulated in detail in Section \ref{sec2}, 
in particular focusing on the mixity of delamination modes.
We confine ourselves to {\it quasistatic} problems 
(i.e.\ inertia neglected) at {\it small strains} and, just for the sake
of notational simplicity, we restrict the analysis to the case of
two (instead of several) visco-elastic bodies $\Omegaone$ and $\Omegatwo$ 
glued together along the {\it contact interface} $\GC$. 
We assume an elastic response of the adhesive, and then one speaks
about {\it adhesive contact}. The elastic response in the adhesive is 
assumed
linear, being determined by the (positive-definite) matrix of elastic
moduli $\bbA$, and the adhesive layer itself is assumed infinitesimally thin. 
At a current time instant, the ``volume fraction'' of
debonded molecular links will be ``macroscopically'' described by
the {\it scalar delamination parameter} $z:\GC\to[0,1]$, which  can
be referred to the modeling approach by M.\,Fr\'emond, see \cite{Fre82,Fre87}.
The state $z(x)=1$ means that the adhesive at $x\in\GC$  is
still $100\%$ undestroyed and thus fully effective, while the
intermediate state $0<z(x)<1$ means that the fraction $1-z(x)$ of  molecular
links  have already been broken but the remaining portion $z(x)$  is still effective,
and eventually $z(x)=0$ means that the interface is already completely
debonded at $x\in\GC$.
As already used e.g.\ in~\cite{KoMiRo},
in some simplification, 
it is assumed that a specific phenomenologically
prescribed energy $a$ (in J/m$^2$, in 3-dimensional situations) is
needed to break the macromolecular structure of the adhesive,
independently of the rate of this process. Thus, delamination is a
\emph{rate-independent} and activated phenomenon, ruled  by the
maximum dissipation principle, and we shall therefore consider a
rate-independent flow rule for $z$.
We will consider the whole energy spent for the delamination as dissipated;
for a more general model cf.\ Remark~\ref{rem-a0-a1} below.

Let us now emphasize main new features of the model, i.e.\ its mixity-sensitivity.  
A standard engineering approach as e.g.~in
\cite{HutSuo92MMCL,TMGCP10ACTA,TMGP11BEMA} is to make the activation
energy $a=a(\psiG)$ depend on the so-called \emph{mode-mixity angle} $\psiG$. 
For instance, if $\norm=(0,0,1)$ at some $x \in \GC$ (with $\norm$ the unit 
normal to $\GC$, oriented from $\Omegaone$ to $\Omegatwo$), and 
$\bbA={\rm diag}(\kappa_{\rm n},\kappa_{\rm t},\kappa_{\rm t})$, the
mode-mixity angle is defined as
$\psiG=\psiG(\JUMP{u}):={\rm arc\,tan}(k^{1/2}|\JUMP{u}_{\rm t}|/|\JUMP{u}_{\rm n}|)$
with $k=\kappa_{\rm t}/\kappa_{\rm n}$
where $\JUMP{u}_{\rm t}$ and $\JUMP{u}_{\rm n}$ stand for
the tangential and the normal traction; i.e.\ 
the jump of displacement across the boundary $\GC$ decomposes 
as $\JUMP{u}= \JUMP{u}_{\rm n}\norm+\JUMP{u}_{\rm t}$, with
$\JUMP{u}_{\rm n}=\JUMP{u}\cdot\norm$.
In fact, to avoid discontinuity of such formula at $0$,
rather a suitable regularization of 
this mode-mixity angle should be taken, e.g.\
\begin{align}\label{a-mixity-dependence-}
\psiG(\JUMP{u})={\rm arc\,tan} \sqrt{\frac{\kappa_{\rm
t}|\JUMP{u}_{\rm t}|^2} {\kappa_{\rm n}|\JUMP{u}_{\rm n}|^2+\epsilon}}
\quad \text{with a small $\epsilon>0$.}
\end{align}
The coefficient  $\kappa_{\rm t}$ is often smaller than $\kappa_{\rm
n}$, and a typical phenomenological form of $a$ used in
engineering \cite{HutSuo92MMCL} is, e.g.,
\begin{align}\label{a-mixity-dependence}
a(\psiG):=a_{\rm I}\big(1+{\rm tan}^2((1{-}\lambda)\psi_{\rm
G})\big);
\end{align}
cf.\ also \cite{BanAsk00NFCI} for a similar formula.
In \eqref{a-mixity-dependence},  $a_{\rm I}=a(0)$ is the activation
threshold for the delamination mode I and $\lambda$ is the so-called
\emph{delamination-mode-sensitivity} parameter.  Note that  a
moderately strong
delamination-mode sensitivity occurs when the ratio $a_{\rm
II}/a_{\rm I}$ is about 5-10 where $a_{\rm II}=a(90^\circ)$ is the
activation threshold for the pure
delamination mode II. Then, one has $\lambda$ about 0.2-0.3; cf.~\cite{TMGP11BEMA}.

Mathematical troubles arising in the analysis of the 
system of partial differential inequalities
for mixity-sensitive delamination model are caused 
by an inelastic rate-independent process on the boundary $\GC$ along which delamination performs.
%
Mixity-dependence of the dissipation makes  the model  {\it
non-associative}, in contrast to the mixity-insensitive case and to
another model involving a interfacial plasticity as an additional internal 
variable, recently devised and analyzed in an isothermal case in
\cite{tr+mk+jz,RoMaPa13QMMD,RoPaMa??LSAQ}. In this work, quadratic viscosity 
dissipation energy coupled with a stored energy density, asymptotically growing 
faster than the spatial dimension and exhibiting separate convexity 
in strain and the delamination variable allow us to prove 
the existence of a suitable weak solution. The viscosity has 
thus  less growth than the stored energy, which, unfortunately,
does not seem to allow for the mechanical  energy 
conservation, cf.\ Remark~\ref{rem-engr} below.



The plan of the paper is as follows: in Section~\ref{sec2} we
formulate the initial-boundary-value problem in its classical
formulation, following essentially 
\cite{rr+tr,rr+tr2}.
After a suitable weak formulation based on the concept of
local solutions due to \cite{Miel11altDEMF,ToaZan09AVAQ}, in Sec.~\ref{sec-ent-trans}
we device a semi-implicit discretization in time combined with
a finite-element discretization in space and prove its numerical
stability, i.e.\ suitable a-priori estimates. The efficient numerical 
implementation (based on recursive quadratic programming combined
with elimination of bulk nodes by a boundary-element method)  
as well as illustrative 2-dimensional computational simulations are performed 
in Sec.~\ref{sec-simul}. Eventually, the convergence of the approximate solutions 
towards weak solutions is outlined in Sect.~\ref{sec-main-res} by using the
already derived a-priori estimates. In particular, we prove unconditional convergence of 
discrete solutions to a weak solution to the model whenever the mesh size and the 
time step tend to zero. 


Let us emphasize that, to our best knowledge, this article represents a first attempt 
to pose the standard engineering model for the mixity-sensitive delamination
of simple materials (together with a physically relevant concept of its 
solution) and to devise an efficiently implementable algorithm in a way which 
simultaneously allows for a rigorous mathematical support as far as numerical
stability and guaranteed convergence.

\section{The model in its classical formulation}\label{sec2}
Hereafter, we suppose that the visco-elastic/inelatic-adhesive structure occupies 
a bounded Lipschitz domain $\Omega\subset\R^d$ composed from (for notational simplicity
only) two visco-elastic bodies, denoted by $\Omegaone$ and $\Omegatwo$, 
glued together on a common conctact boundary, denoted by $\GC$, which
represents a prescribed
interface where delamination may occur. 
This means we consider
\[
\Omega = \Omegaone \cup \GC \cup \Omegatwo\,,
\]
with $\Omegaone$ and $\Omegatwo$
disjoint Lipschitz subdomains. 
We  denote by $\vec{n}$ the
outward unit normal to $\partial \Omega$, and by $\norm$
 the unit normal to $\GC$, which we consider
oriented from $\Omegaone$ to $\Omegatwo$. Moreover, given  $v$
defined on $\Omega {\setminus} \GC$,
$v^+$ (respectively, $v^-$)
signifies the restriction of $v$ to $\Omegaone$ (to $\Omegatwo$,
resp.).
 We further suppose  that the boundary of $\Omega$ splits as
\[
\partial \Omega = \GD\cup \GN\,,
\]
with $\GD$ and $\GN $ open subsets in the relative topology of
$\partial\Omega$,  disjoint one from each other and each of them
with a smooth (one-dimensional) boundary. Considering $T>0$ a fixed
time horizon, we  set
\begin{displaymath}
Q:=(0,T){\times}\Omega, \quad \Sigma:= (0,T){\times}\partial
\Omega, \quad \SC\!:= (0,T){\times}\GC , \quad \Sdir\!:= (0,T){\times}
\GD, \quad \Snew\!:= (0,T){\times}\GN.
\end{displaymath}
For readers' convenience, let us summarize the basic notation used
in what follows:

\vspace{.7em}

\hspace*{-1.4em}\fbox{
\begin{minipage}[t]{0.46\linewidth}
\small

$d=2,3$ dimension of the problem,

$u:\Omega{\setminus}\GC \to\R^d$ displacement,


$z:\GC\to[0,1]$ delamination variable,

$e=e(u)=\frac12\nabla u^\top\!+\frac12\nabla u$ small-strain tensor,

$\JUMP{u}= u^+|_{\GC} - u^-|_{\GC}\ $ jump of $u$ across $\GC$,

$\sigma$  stress tensor,

$\psiG$ mode-mixity angle,






\end{minipage}\ \
\begin{minipage}[t]{0.49\linewidth}\small

$\bbC:\R_{\mathrm{sym}}^{d\times d}\to\R_{\mathrm{sym}}^{d\times d}$
nonlinear elastic Hook law,

$\bbD\in\R^{d^4}$ viscosity constants,


$\bbA\in\R^{d\times d}$
elastic coefficients of the adhesive,



$\alpha=\alpha(\JUMP{u})$ energy (per area) dissipated on $\GC$,




$\FRM:Q\to\R^d$  applied bulk force,

$w_{\rm D}
$  prescribed
boundary displacement,

$\fRM:\Snew\to\R^d$  applied traction.






\end{minipage}\medskip
}

\vspace{-.5em}

\begin{center}
{\small\sl Table\,1.\ }
{\small
Summary of the basic notation used thorough the paper. 
}
\end{center}

\noindent The {\it state} is formed by the couple $(u,z)$.
We use {\it Kelvin-Voigt's rheology} and, rather for mathematical
reasons to facilitate analysis in multidimensional cases,
a (possibly only slightly) nonlinear static response.
Hence we assume the \emph{stress} 
$\sigma:(0,T)\times\Omega\rightarrow\R^{d\times d}$ in the form:
\begin{align}
\sigma=\sigma(u,\DT{u}) :=\!\!\!\ddd{\bbD
e(\DT{u})}{viscous}{stress}\!\!\!
+\!\!\!\ddd{\bbC(e(u))}{elastic}{stress}\!\!\!,
\end{align}
Furthermore, we shall denote by $T=T(u,v)$  the traction
stress on some two-dimensional surface $\Gamma$
(later, we shall take either $\Gamma=\GC$ or $\Gamma=\GN$), i.e.
\begin{align}\label{T-stress}
T(u,\DT{u}):=\sigma(u,\DT{u})\big|_{\Gamma}\vec{n}\,,
\end{align}
where of course we take as $\vec{n}$ the unit normal  $\norm$ to $\GC$,
if $\Gamma=\GC$. Its normal and tangential components are defined on $\GC\cup\GN$ 
respectively by the formulas 
\begin{align}\label{def-of-T}
T_{\rm n}(u,\DT{u})=\vec{n} \cdot\sigma(u,\DT{u})\big|_{\Gamma}\vec{n}
\qquad\text{and}\qquad
T_{\rm t}(u,\DT{u})=T(u,\DT{u})-T_{\rm n}(u,\DT{u})\vec{n}.
\end{align}

We address 
the standard frictionless Signorini
conditions on $\GC$ for the displacement $u$.  
\paragraph{Classical formulation of the adhesive contact problem.}
 Beside the force equilibrium coupled with the heat equation inside
   $\Q{\setminus}\SC$ and supplemented  with standard boundary
  conditions, we have two complementarity problems on $\SC$. Altogether, we
  have the boundary-value problem
\begin{subequations}\label{adhes-classic}
\begin{align}
\label{eq6:adhes-class-form1} & 
\mathrm{div}\big(\bbD e(\DT{u})+\bbC(e(u))
\big)+\FRM=0, 
&&\text{in }Q{\setminus}\SC,
\\
\label{eq6:adhes-class-form2} &u=0 &&\text{on }\Sdir,\hspace{1.2em}
\\\label{eq6:adhes-class-form3-bis}
&T(u,\DT{u})
=\fRM &&\text{on }\Snew,\hspace{1.2em}
\\\label{adhes-form-d1}
& \JUMP{\bbD e(\DT{u})+ \bbC(e(u))
}\norm 
=0
 &&\text{on
}\SC,\hspace{1.2em}
\\
\label{adhes-form-d2-} &
T_{\rm t}(u,\DT{u})
+z\big(\mathbb A u{-}\big((\mathbb A u){\cdot}\norm\big)\norm\big)=0
&&\text{on }\SC,\hspace{-.15em}
\\
\label{adhes-form-d2} &\JUMP{u}{\cdot}\norm\ge0  
\ \ \ \ \text{ and }\ \ \ \ T_{\rm n}(u,\DT{u}) 
+z(\bbA\JUMP{u}){\cdot}\norm
\ge0
&&\text{on }\SC,\hspace{1.2em}
\\
\label{adhes-form-d4} &
\big(T_{\rm n}(u,\DT{u})
{+}z(\bbA\JUMP{u})\norm
\big)(\JUMP{u}{\cdot}\norm)=0
&&\text{on }\SC,\hspace{1.2em}
\\\label{adhes-form-d6}
&\DT{z}\le0
\ \ \ \ \text{ and }\ \ \ \ d\le 
\alpha(\JUMP{u}) \ \ \ \ \text{ and }\ \ \ \ 
\DT{z} \left(d - \alpha(\JUMP{u})\right) =0
&&\text{on }\SC,\hspace{1.2em}
\\\label{adhes-form-d8}
& d\in \partial I_{[0,1]}(z)+
\mbox{$\frac12$}\bbA\JUMP{u}{\cdot}\JUMP{u} &&\text{on
}\SC.\hspace{1.2em}
\end{align}
\end{subequations}
%
As to the involved  symbols,
we assume that
\begin{subequations}\label{ass-on-K-C-D-H-G}
\begin{align}
\label{posit}&
\varphi:\R_{\mathrm{sym}}^{d\times d}\to\R\text{ convex smooth};\,\exists\,\varepsilon_0,\varepsilon_1>0, p>d:\,\varepsilon_1(1+|e|^p)\ge\varphi(e)\ge\varepsilon_0(|e|^p\!-1),\\
\label{posit-C}& 
\bbC(e):=\varphi'(e):=\frac{\partial\varphi(e)}{\partial e}\quad\text{ elastic stress tensor at strain $e$},
\\\label{posit-D}&
\bbD:\R_{\mathrm{sym}}^{d\times d}\to
\R_{\mathrm{sym}}^{d\times d}\quad\text{ linear 
positive definite,}
\\\label{posit-A}&
\bbA:\R^d\to
\R^d\quad\text{ linear 
positive semidefinite}.
\end{align}\end{subequations}

The complementarity problem \eqref{adhes-form-d2}--\eqref{adhes-form-d4} 
describes the Signorini unilateral contact. 
%
The complementarity problem~\eqref{adhes-form-d6}--\eqref{adhes-form-d8}
corresponds to the \emph{flow rule} governing the evolution of $z$:
\begin{equation}
\label{e:reactivation-2}
\partial \ind_{(-\infty,0]}(\DT{z}) + \partial\ind_{[0,1]}(z) +
\mbox{$\frac12$}\bbA\JUMP{u}{\cdot}\JUMP{u}
\ni\alpha(\JUMP{u})\qquad\text{in $\SC$,}
\end{equation}
with the indicator functions $\ind_{(-\infty,0]},\, \ind_{[0,1]}: \R \to
[0,+\infty]$ and their (convex analysis) subdifferentials 
$\partial\ind_{(-\infty,0]},\,\partial\ind_{[0,1]}: \R \rightrightarrows\R$. 
The energetics of the model is formally:
\begin{align}\nonumber
  &\frac{\d}{\d t}\bigg(\!\!\ddd{\int_{\Omega{\setminus}\GC}\!
\varphi(e(u)) 
\,\d x}{
elastic
energy}{in the bulk}\!\!\!\!\! 
+ \!\!\!\!\!\ddd{\int_{\GC}\!\! \frac z2{\dela}
\JUMP{u}{\cdot}\JUMP{u}
\,\d S}{elastic energy}{in the adhesive}\!\!\bigg)
+\!\!\!\!\!\ddd{\int_{\Omega{\setminus}\GC}\!\frac12\bbD e(\DT u){:}e(\DT u)
\,\d x}{rate of viscous}{dissipation in the bulk}\!\!\!\!\!\!
\\&\label{total-energy-}
\hspace{12em} + \!\!\!\!\!\!\!\!\!
\ddd{\int_{\GC}\!\!\alpha(\JUMP{u})\DT z\,\d S}
{rate of dissipation by dela-}{mination of the adhesive}\!\!\!\!
   =\!\!\!\!\ddd{\int_{\Omega}\!
\FRM{\cdot}\DT{u}\,\d x}
{power of bulk}{mechanical load} \!\!\!\!
+\!\!\!\ddd{\int_{\GN}\!\!\fRM{\cdot}\DT{u}\,\d S}
{power of surface}{mechanical load}\!\!\!.
\end{align}
For more details about derivation of the model we refer to \cite{rr+tr,rr+tr2}.
We will consider the initial-value problem for \eqref{adhes-classic} by
prescribing the initial condition
\begin{align}\label{IC}
u(0)=u_0 \quad \aein \ \Omega, \qquad z(0)=z_0 \quad
\aein \ \GC.
\end{align}

\begin{remark}[{\it Stored energy increased by delamination}]\label{rem-a0-a1}
\upshape
The energy needed for the delamination can be alternatively understood 
as contributing to the stored energy. This reflects the fact that
any new surface represents some additional stored energy. 
In the isothermal unidirectional delamination, this alternative
concept is mechanically equivalent. Yet, if temperature variations
are considered, then it makes a difference because the stored energy
variation does not contribute to the heat production. In fact,
rather both parts (i.e.\ dissipative and stored) of the energy spent
for delamination should more realistically be considered, cf.\ 
\cite{rr+tr2}. 
Also, if 
a bi-directional evolution of delamination (involving healing)
would be considered, then the contribution to the stored energy
becomes especially important because it just facilitates the driving
force for possible healing, cf.\ \cite{tr-os-rv}.
\end{remark} 

\begin{remark}[{\it Dynamical problems}]\label{rem-inertia}
\upshape
In some applications/regimes inertial forces cannot be neglected and then
\eqref{eq6:adhes-class-form1} takes the form
\begin{align}
\varrho\DDT{u}-\mathrm{div}\big(\bbD e(\DT{u})+\bbC
e(u)-\mathrm{div}\,\mathfrak h\big)=\FRM
\end{align} 
with $\varrho>0$ mass density. Implicit discretization of this
term is relatively easy to be incorporated in the analysis
if 
a generalized concept of solution without energy preservation is accepted, 
cf.\ \cite{rr+tr2}. 
Here e.g.\ \eqref{semi-impl+1} augments by the term 
$\tau^{-2}\int_\Omega\frac12\varrho|u{-}2u_\tau^{k-1}{+}u_\tau^{k-2}|^2\d x$.
Yet, it is well known that the implicit discretization of the inertial
term is unsuitable for computational simulations due to spurious 
numerical attenuation and efficient calculations of wave propagation 
needs more sophisticated formulas. On the other hand, leaving the 
energy preservation out, we can also afford $\bbD=0$ because the 
inertial term controls $\JUMP{u}$ ``compactly'' in $C(\overlineSC)$
via Aubin's-Lions' theorem, thus we get hyperbolic inviscid delamination problem.
\end{remark} 


\begin{remark}[{\it Cohesive contacts}]\label{rem-inertia-1}
\upshape
We can also consider $z\dela{+}z^2\bbB$ instead of $z\dela$ in (\ref{adhes-classic}g-i) 
and $(\frac12\dela{+}z\bbB)\JUMP{u}{\cdot}\JUMP{u}-\kappa\Delta z$
instead of $\frac12\dela\JUMP{u}{\cdot}\JUMP{u}$ in \eqref{adhes-form-d8},
which would be based on the stored energy with the boundary term of the type 
$$
\int_{\GC}\frac12(z\dela{+}z^2\bbB)
\JUMP{u}{\cdot}\JUMP{u}+\kappa|\nabla_{_{\rm S}}z|^2\,\d S
$$ 
with $\bbB$ positive semidefinite and $\kappa>0$;
for a more detailed discussion about this quadratic cohesion model 
cf.\ \cite[Sect.6.1]{tr+mk+jz} and for analysis cf.\ \cite{BBR1}.
Here, it leads to two quadratic mathematical programs after the semi-implicit 
discretization and numerical analysis works for P1-discretization of $z$ and
2-dimensional problems simply by a mutual recovery sequence
$\widetilde z_k:=(\widetilde z-\|z_k{-}z\|_{C(\overlineGC)})^+$ suitably adjusted 
to spatial discretization, cf.\ the proof of stability of the limit 
below (while for 3-dimensional problems more sophisticated 
damage-type construction by M.Thomas et al.\ \cite{Thom10PhD,ThoMie10DNEM} 
would be needed).
\end{remark} 

\section{Weak formulation and semi-implicit discretization}\label{sec-ent-trans}

%
We will use the standard notation $W^{1,p}(\Omega)$ for the Sobolev space 
of functions having the derivatives in the Lebesgue space $L^p(\Omega)$. 
If valued in $\R^n$ with $n\ge2$, we will write $W^{1,p}(\Omega;\R^n)$, and 
furthermore, if $p=2$, we 
use the shorthand notation $H^1(\Omega;\R^n)=W^{1,2}(\Omega;\R^n)$.
Moreover, we will adopt the notation
\[
\begin{aligned}\label{def-of-W}
&  W_{\GD}^{1,p}(\Omega {\setminus} \GC ;\R^d):=\big\{ \testu
\in W^{1,p}(\Omega {\setminus} \GC ;\R^d): \ \ \testu =0 \ \
\text{on $\GD$}\big\}\,.
\end{aligned}
\]
For $X$ a (separable) Banach space, we denote by 
$C_{\rm w}([0,T];X)$ and $BV([0,T];X)$ the Banach spaces 
of weakly continuous functions $[0,T]\to X$ and of the functions that have 
bounded variation on $[0,T]$, respectively. Notice that these functions 
are defined everywhere on $[0,T]$.

Hereafter, the external mechanical
loading  $\FRM$ and $\fRM$
will be qualified 
\begin{subequations}\label{hypo-data}
\begin{align}
  \label{eFFe1}
&\FRM \in\begin{cases} L^2 (0,T; L^{6/5}(\Omega; \R^d)) & \mbox{ if $d=3$},\\
 L^2(0,T;L^q(\Omega;\R^d))\, ,\, q>1 & \mbox{ if $d=2$};
 \end{cases}
\\
& \label{eFFe2}\fRM \in  \begin{cases} L^2(0,T; L^{4/3}(\GN;\R^d)) & \mbox{ if $d=3$},\\
L^2(0,T; L^{q}(\GN;\R^d)) \, ,\, q>1 & \mbox{ if $d=2$};
\end{cases}
\end{align}
\end{subequations}
cf.\ also Remark~\ref{rem-F} below.
As for the initial data, we impose the following
\begin{subequations}
\label{hyp-init}
\begin{align}
& \label{uzero} u_0 \in W_{\GD}^{1,p}(\Omega{\setminus}
\GC;\R^d)\,, \quad \JUMP{u_0}{\cdot}\norm \ge 0 \ \ \text{on $\GC$,}
\\
& \label{zzero} z_0 \in L^\infty(\GC), \qquad 0 \leq z_0 \leq 1 \ \
\text{a.e. on}\, \GC\,.
\end{align}
\end{subequations}
We will use the abbreviation for the stored energy $\Phi$
and  the  dissipation rate  $\calD$:
\begin{align}
&\Phi(u,z):=\begin{cases}
\displaystyle{
\int_{\Omega{\setminus}\GC}\!\!\!
\varphi(e(u))
\,\d x
+\int_{\GC}\!
\frac12z\dela\JUMP{u}{\cdot}\JUMP{u}
\,\d S} & \text{if $\JUMP{u}{\cdot}\norm\ge0$ 
and $0{\le}z{\le}1$ on $\GC$,}
\\[-.2em] + \infty & \text{otherwise, and}\end{cases}
\label{8-1-k}
\\&
\label{DISS} \calD\big(u;\DT u,\DT z):=
\begin{cases}
\displaystyle{
\int_\Omega\bbD e(\DT{u}){:}e(\DT{u})\,\d x
+\int_{\overlineGC}\!\alpha(\JUMP{u})|\DT z|\,\d S}  & \text{if
$\DT z\le0$ a.e. in $\GC$,}
\\[-.2em] + \infty & \text{otherwise.}
\end{cases}
\end{align}

\begin{definition}[{\it Weak solution}]\label{def4}
\upshape
 Given an initial data
$(u_0,z_0)$ satisfying \eqref{hyp-init},
 we call a couple $(u,z)$ a weak solution to the
Cauchy problem for system~\eqref{adhes-classic} if\nolinebreak
\begin{subequations}
\label{reguu}
\begin{align}
\label{reguu1} & \!\!\! \!\!\!u \in
C_{\rm w}([0,T];W_{\GD}^{1,p}(\Omega{{\setminus}}\GC;\R^d))
\cap H^1(0,T;H^1(\Omega{{\setminus}}\GC;\R^d)),
\\
\label{reguz} &  \!\!\! \!\!\!z \in L^\infty (\SC)\, \cap\,
BV([0,T];L^1(\GC))\,,\  \text{
$z(\cdot,x)$ nonincreasing on $[0,T]$ for a.a. } x\!\in\!\GC, 
\end{align}
\end{subequations}
and the couple $(u,z)$  complies, besides the initial condition \eqref{IC}, 
with:
\\
 \ITEM{(i)}{(weak formulation of the) momentum inclusion, i.e.:}\nolinebreak
\begin{subequations}\label{weak-delam}\begin{align}
\label{constraints-delam}
& \JUMP{u}{\cdot}\norm\ge0\ \ \text{on $\SC$,} \quad \text{ and}
\\
\nonumber
 &
\int_{Q\setminus\SC}\!\!\!\!
\big(\bbD e(\DT{u}){+}\bbC(e(u))\big){:}e(\testu{-}u)
\,\d x\d t
+\int_{\SC}\!\!\!
z\bbA\JUMP{u}
{\cdot}\JUMP{\testu{-}u}\d S\d t
\\&\qquad\qquad\qquad\qquad\qquad
\ge
\int_Q\!\FRM{\cdot}(\testu{-}u)\,\d x\d t
+\int_{\Snew}\!\!\fRM{\cdot}(\testu{-}u)\,\d S\d t
\label{e:weak-momentum-var}
\end{align}
\ITEM{}{for all $\testu$ in
$L^2(0,T;W_{\GD}^{1,p}(\Omega{\setminus}\GC;\R^d))$ with
$\JUMP{\testu}{\cdot}\norm\ge0$  on~$\SC$,
}
\ITEM{(ii)}{energy inequality for almost all time instant  $t_1<t_2$, $[t_1,t_2]\subset [0,T]$:}
\begin{align}
\Phi\big(u(t_2),z(t_2)\big)&+\int_{t_1}^{t_2}\!\calD(u;\DT u,\DT z)\,\d t\nonumber\\[-.3em]
&\le
\Phi\big(u(t_1),z(t_1)\big)
+\int_{t_1}^{t_2}\!\int_\Omega\!\FRM{\cdot}\DT u\,\d x\d t
+\int_{t_1}^{t_2}\!\int_{\GN}\!\!\fRM{\cdot}\DT u\,\d S \d t,
 \label{total-energy}
\end{align}
 \ITEM{(iii)}{
semistability for a.a.~$t\in (0,T)$}
\begin{align}
 \label{semistab} &\forall \tilde{z}\in
L^\infty(\GC):\qquad \Phi\big(u(t),z(t)\big)\le\Phi\big(u(t),\tilde
z\big) +\calD\big(u(t);0,\tilde z-z(t)\big).
\end{align}
\end{subequations}
\end{definition}

\bigskip

\begin{remark}\upshape 
Due to cancellation of bulk terms, \eqref{semistab} means just
$\int_{\GC}(z(t){-}\tilde z)(\dela\JUMP{u(t)}{\cdot}\JUMP{u(t)}
-2\alpha(\JUMP{u(t)}))\,\d S\le0$ for all 
$\tilde z\!\in\! L^\infty(\GC)$ such that $0\!\le\!\tilde z\!\le\!z(t)$,
which can be disintegrated so that \eqref{semistab} is equivalent to
\begin{align}\label{semi-stabil+}
z(t,x)\dela\JUMP{u(t,x)}{\cdot}\JUMP{u(t,x)}\le2\alpha(\JUMP{u(t,x)})
\ \ \text{ or }\ \ z(t,x)=0\qquad \text{ for a.a.\ $x\in\GC$.}
\end{align}
In 
\cite{Miel05ERIS,MieThe04RIHM}, 
a global stability
condition combined with energy conservation was shown to provide the
correct ``weak'' formulation of rate-independent flow rules. Due
to the viscosity in the bulk, here the semistability
\eqref{semistab} plays the role of the global stability condition 
of 
\cite{Miel05ERIS,MieThe04RIHM}. 
Moreover, 
here we do not require the energy conservation \eqref{total-energy-}
but we only assume energy inequality \eqref{total-energy} between 
varying time instances $t_1<t_2$, which is the general concept of 
so-called {\it local solutions} invented
for purely rate-independent systems for a special crack problem in 
\cite{ToaZan09AVAQ} and further generally investigated in \cite{Miel11altDEMF},
and proved to coincide with the concept of conventional weak solutions 
in \cite[Prop.2.3]{Roub??MDLS}. Although this concept is very wide in general, 
here the viscosity together with the convexity of the stored energy in terms of $z$ 
ensures good selectivity of this concept, cf.\ \cite[Prop.5.2]{tr2}.
This viscosity/convexity attribute also excludes the undesired effect of too-early 
delamination unphysically sliding to less dissipative Mode I,
which may occur in purely elastic model if energy conservation 
would be forced, cf.\ \cite{RoMaPa13QMMD} in contrast to 
\cite{RoPaMa??LSAQ}.
Also we point out that, disregarding the only one-sided inequality, 
\eqref{total-energy} is the integrated version of the total 
energy balance~\eqref{total-energy-}. 
%
%


\end{remark}


To solve the the initial-boundary value problem 
\eqref{adhes-classic} and \eqref{IC} numerically, we 
must make some discretization both in time and in space. 
Rather as an example, let us consider P1-elements for $u$ and 
P0-elements for $z$. Assuming polygonal domains $\Omegaone$ and 
$\Omegatwo$, we use 
a spatial discretization by considering a triangulation $\mathscr{T}_h$
of $\Omega{\setminus}\GC$ with a mesh size $h>0$ and define the 
finite-dimensional subspaces 
\begin{subequations}\label{Uh-Vh}\begin{align}\label{Uh}
&V_h:=\big\{v\in W^{1,\infty}(\Omega;\R^d);\ \ \forall S\in\mathscr{T}_h:
\ \ v|_S\text{ affine}\big\},
\\\label{Zh}
&Z_h:=\big\{z\in L^\infty(\GC);\ \ \forall S\in\mathscr{T}_h:
\ \ z|_{\bar S\cap\GC}\text{ constant}\big\}.
\end{align}\end{subequations}
Moreover, we make the {\it time-discretization} by using 
a suitable {\it semi-implicit scheme} using a popular
fractional-step-like strategy, cf.\ also \cite[Remark~8.25]{NPDE_roubicek}.
In contrast to anisothermal situation in \cite{rr+tr2}, here this leads to
alternating variational problems which are even convex, which allows for
a constructive solution.
Using an equidistant partition of the time interval $[0,T]$ with a time 
step $\tau>0$, 
we seek $u_{\tau h}^k\in V_h$ and $z_{\tau h}^k\in Z_h$
such that $\JUMP{u_{\tau h}^k}{\cdot}\norm\ge0$ on $\SC$ and
\begin{subequations}\label{disc}
\begin{align}\nonumber
 &\int_{\Omega\setminus\GC}\!\!\Big(\bbD e\big(\frac{u_{\tau h}^k{-}u_{\tau h}^{k-1}}\tau\big)
{+}\bbC(e(u_{\tau h}^k))\Big){:}e(\testu{-}u_{\tau h}^k)\,\d x
+\int_{\GC}\!\!\!
z_{\tau h}^{k-1}\bbA\JUMP{u_{\tau h}^k}
{\cdot}\JUMP{\testu{-}u_{\tau h}^k}\d S
\\\label{disc-u}
&\hspace{12em}\ge\int_\Omega\!\FRM_\tau^k{\cdot}(\testu{-}u_{\tau h}^k)\,\d x
+\int_{\GN}\!\!\fRM_\tau^k{\cdot}(\testu{-}u_{\tau h}^k)\,\d S
\\&\nonumber
\Phi\big(u_{\tau h}^k,z_{\tau h}^k\big)
+\calD\Big(u_{\tau h}^k;\frac{u_{\tau h}^k{-}u_{\tau h}^{k-1}}\tau,
\frac{z_{\tau h}^k{-}z_{\tau h}^{k-1}}\tau\Big)\le\Phi\big(u_{\tau h}^{k-1},z_{\tau h}^{k-1}\big)
\\\label{disc-z-1}
&\hspace{12em}
+\int_\Omega \FRM_\tau^k{\cdot}\frac{u_{\tau h}^k{-}u_{\tau h}^{k-1}}\tau\,\d x 
+\int_{\GN}\!\!\fRM_\tau^k{\cdot}\frac{u_{\tau h}^k{-}u_{\tau h}^{k-1}}\tau\,\d S,
\\&\label{disc-z-2}
\forall \tilde{z}\in Z_h:\qquad \Phi\big(u_{\tau h}^k,z_{\tau h}^k\big)
\le\Phi\big(u_{\tau h}^k,\tilde z\big) +\calD\big(u_{\tau h}^k;0,\tilde z-z_{\tau h}^k\big).
\end{align}
\end{subequations}
with 
$\FRM_\tau^k=\tau^{-1}\int_{(k-1)\tau}^{k\tau}\FRM(s)\,\d s$, $\fRM_\tau^k=\tau^{-1}\int_{(k-1)\tau}^{k\tau}\fRM(s)\,\d s$,
and proceeding recursively for $k=1,... T/\tau\in\N$ with 
starting for $k=1$ from 
\begin{align}
u_{\tau h}^0=u_0\qquad\text{ and }\qquad z_{\tau h}^0=z_0.
\end{align}
The adjective ``semi-implicit'' is related with usage of $z_{\tau h}^{k-1}$ 
in (\ref{disc-u}) instead of $z_{\tau h}^k$ which would lead to a fully implicit 
formula.
Such usage  of $z_{\tau h}^{k-1}$ leads to the decoupling of the problem:
first we can solve (\ref{disc-u}) for $u_{\tau h}^k$ and, 
only after it, the rest of (\ref{disc}b,c) for $z_{\tau h}^k$. 
Note that, in a simple way, we discretized rather the weak formulation 
\eqref{weak-delam} than the classical formulation 
\eqref{adhes-classic} where we would have faced technical problems 
e.g.\ with the interaction of a piece-quadratic 
$\bbA\JUMP{u_{\tau h}^k}{\cdot}\JUMP{u_{\tau h}^k}$ with piecewise constant 
$z_{\tau h}^k{-}z_{\tau h}^{k-1}$ and a general nonlinear $\alpha(\JUMP{u_{\tau h}^k})$
in (\ref{adhes-classic}h-i).

On top of it, we can employ the variational structure of both 
decoupled problems, cf.\ also \cite[Remark 8.25]{NPDE_roubicek}. 
To cope with the constraints more explicitly, we introduce the smooth
stored energy $\calE$ and the 
dissipation (pseudo)potentials $\calR_1$ and $\calR_2$ defined here by
\begin{subequations}\label{E-R}\begin{align}\label{M}
&\calE(t,u,z)=\int_{\Omega}\!
\varphi(e(u))
-\FRM(t){\cdot}u\,\dd x
+\int_{\GC}\!\frac12z\mathbb A\JUMP{u}{\cdot}\JUMP{u}\,\d S
-\int_{\GN}\!\!\fRM(t){\cdot}u\,\d S,
\\&\calR_1(u;\DT z)=-\int_{\GC}\!\!\alpha(\JUMP{u})\DT z\,\d S\ ,\qquad\quad
\calR_2(\DT u)=\int_\Omega\frac12\bbD e(\DT u){:} e(\DT u)\,\d x.
\end{align}\end{subequations}
Note that the constraints $\JUMP{u}{\cdot}\norm\ge0$, $0\le z\le 1$, 
and $\DT z\le0$, originally contained in the stored energy $\Phi$ and 
the dissipation rate $\calR$ in
\eqref{8-1-k} and \eqref{DISS}, are now included in 
\eqref{semi-impl+} below so that
we can equivalently use the smooth functionals $\calE(t,\cdot,\cdot)$
and $\calR_1(u;\cdot)$ and $\calR_2$ in \eqref{E-R}. Also 
note that $\calR(u;\DT u,\DT z)=\calR_1(u;\DT z)+2\calR_2(\DT u)$
and $\calR_1(u;\cdot)$ is  degree-1 homogeneous
so that the factor $\tau$ does not show up in the functional 
in \eqref{semi-impl+2}, in contrast to the degree-2 homogeneous functional
$\calR_2(\cdot)$ in \eqref{semi-impl+1}.
We thus obtain two convex minimization problems: first, we are to solve
\begin{subequations}\label{semi-impl+}\begin{align}\label{semi-impl+1}
&\left.\begin{array}{ll}
\text{minimize}&\displaystyle{\calE(k\tau,u,z_{\tau h}^{k-1})
+\tau\calR_2\Big(\frac{u{-}u_{\tau h}^{k-1}}{\tau}\Big)}
\\[.3em]
\text{subject to}&u\in V_h,\ \ 
\JUMP{u}{\cdot}\norm\ge0
\end{array}\hspace{2em}\right\}
\intertext{and, denoting its unique solution by $u_{\tau h}^k$, then we solve}
\label{semi-impl+2}
&\left.\begin{array}{ll}
\text{minimize}&\displaystyle{\calE\big(k\tau,u_{\tau h}^k,z\big)
+\calR_1\big(u_{\tau h}^k;z{-}z_{\tau h}^{k-1}\big)
}
\\[.3em]
\text{subject to}&z\in Z_h,\ \ \ 0\le z\le z_{\tau h}^{k-1}.\end{array}
\hspace{2.7em}\right\}
\end{align}\end{subequations}

For $\tau>0$ fixed, we denote the left-continuous and the right-continuous 
\emph{piecewise constants}, and the \emph{piecewise linear} interpolants 
of the discrete solutions $\{\uk \}_{k=1}^{T/\tau}$ by 
$\pwc u{\tau h}:(0,T)\to W_{\GD}^{1,p}(\Omega{\setminus}\GC;\R^d)$, 
$\upwc u{\tau h}:(0,T)\to W_{\GD}^{1,p}(\Omega{\setminus}\GC;\R^d)$, and
$\pwl u{\tau h}:(0,T)\to W_{\GD}^{1,p}(\Omega{\setminus}\GC;\R^d)$
defined by 
\begin{align}\label{interpolants}
\pwc u{\tau h}(t)=\uk,\ \ \upwc u{\tau h}(t)=\ukm,\ \ \pwl u{\tau h}(t)
=\frac{t-t_\tau^{k-1}}{\tau}\uk+\frac{t_\tau^k-t}{\tau}\ukm\ \ \text{ for $t\in(t_\tau^{k-1},t_\tau^k]$.}
\end{align}
In the same way, we shall denote the interpolants of $\{\zk\}_{k=1}^{T/\tau}$, and of
 $\FRM_k^\tau$, and $\fRM_k^\tau$.


Both for supporting convergence analysis (cf.\ Sect.\,\ref{sec-main-res} below)
and for implementation, the important attribute of the above devised 
discrete scheme is its numerical stability, i.e.\ the numerical results
do not exhibit spurious mesh dependency:

\begin{proposition}[Numerical stability of the discretization]\label{prop:apriori}
Let us assume \eqref{ass-on-K-C-D-H-G}, \eqref{hypo-data}, 
\eqref{hyp-init}, 
$\inf\alpha(\cdot)>0$,
$\mathrm{meas}_{d-1}(\partial\Omegaone\cap\GD)>0$, and 
$\mathrm{meas}_{d-1}(\partial\Omegatwo\cap\GD)>0$,
where $\mathrm{meas}_{d-1}$ denotes the $(d{-}1)$-dimensional 
measure on $\Gamma$.
Then, for all $\tau>0$ and $h>0$ and for some constant $S_0>0$ 
independent of $\tau$ and $h>0$, the approximate solutions $(\pwc u{\tau h},
\pwc z{\tau h}, \pwl u{\tau h}, \pwl z{\tau h})$  satisfy
\begin{subequations}\label{a-priori3}
\begin{align}
& \label{a30} \big\|\pwc
u{\tau h}\big\|_{L^{\infty}(0,T;W_{\GD}^{1,p}(\Omega{\setminus}\GC;\R^d))}\le
S_0\,,
\\
& \label{a31} \big\|\pwl u{\tau h}\big\|_{
H^1(0,T;H^1(\Omega{\setminus}\GC;\R^d))}\le S_0\,,
\\&
\label{a32} \big\|\pwc z{\tau h}\big\|_{L^\infty(\SC)}
\leq S_0\,,
\\
 & \label{a32bis}
 \big\|\pwc z{\tau h}\big\|_{BV([0,T];L^1(\GC))} \leq S_0.
\end{align}
\end{subequations}
\end{proposition}
\noindent{\it Sketch of the proof.} We only sketch the calculations for proving
\eqref{a-priori3}, since the argument closely follows the proof of
\cite[Lemma 7.7]{rr+tr}
or also \cite[Lemma 5.6]{rr+tr2}.

A discrete analog of \eqref{total-energy} 
can be obtained by testing the optimality conditions for \eqref{semi-impl+1}
and \eqref{semi-impl+2} respectively by $u_{\tau h}^k{-}u_{\tau h}^{k-1}$ and 
$z_{\tau h}^k{-}z_{\tau h}^{k-1}$ (which, in fact, means plugging $v=u_{\tau h}^{k-1}$ into 
a discrete version of \eqref{e:weak-momentum-var} for 
the former test), and by adding it, benefiting from the cancellation of the 
terms $\pm\Phi(u_{\tau h}^k,z_{\tau h}^{k-1})$ and by the separate
convexity of $\Phi(\cdot,\cdot)$, i.e.\ both $\Phi(u,\cdot)$
and $\Phi(\cdot,z)$ are convex. This gives the estimate
\begin{align}\nonumber
&\Phi(u_{\tau h}^k,z_{\tau h}^k)
+\tau\sum_{l=1}^k\calR\Big(u_{\tau h}^l;\frac{u_{\tau h}^l{-}u_{\tau h}^{l-1}}\tau,
\frac{z_{\tau h}^l{-}z_{\tau h}^{l-1}}\tau\Big)
\\\nonumber&\le\Phi(u_0,z_0)+
\tau\sum_{l=1}^k
\int_\Omega \FRM_\tau^l{\cdot}\frac{u_{\tau h}^l{-}u_{\tau h}^{l-1}}\tau\,\d x 
+\int_{\GN}\!\!\fRM_\tau^l{\cdot}\frac{u_{\tau h}^l{-}u_{\tau h}^{l-1}}\tau\,\d S
\\\nonumber&\le\Phi(u_0,z_0)+\|\FRM_\tau^l\|_{L^{6/5}(\Omega;\R^d)}
\Big\|\frac{u_{\tau h}^l{-}u_{\tau h}^{l-1}}\tau\Big\|_{L^6(\Omega;\R^d)}
+\|\fRM_\tau^l\|_{L^{4/3}(\GN;\R^d)}
\Big\|\frac{u_{\tau h}^l{-}u_{\tau h}^{l-1}}\tau\Big\|_{L^4(\GN;\R^d)},
\\&\le\Phi(u_0,z_0)+C_\delta\|\FRM_\tau^l\|_{L^{6/5}(\Omega;\R^d)}^2
+C_\delta\|\fRM_\tau^l\|_{L^{4/3}(\GN;\R^d)}^2
+\delta\Big\|\frac{u_{\tau h}^l{-}u_{\tau h}^{l-1}}\tau\Big\|_{H^1(\GN;\R^d)}^2,
\label{engr-ineq-disc}\end{align}
where $\delta>0$ and $C_\delta$ depends, beside $\delta$, also on the
norms of the embedding $H^1(\Omega{\setminus}\GC)\subset L^6(\Omega)$ and 
of the trace operator $H^1(\Omega{\setminus}\GC)\to L^4(\GN)$.
Then we choose $\delta>0$ so small that the last term can be
absorbed in the $\calR$-term by using the assumption \eqref{posit-D}.
Then all the a-priori estimates \eqref{a-priori3} easily follow.
$\hfill\Box$

\begin{remark}\label{rem-F}
\upshape 
The integrability in \eqref{hypo-data} designed rather
for $d=3$ can be improved for $d=2$. One can also consider the 
alternative qualification, e.g.\ for $p>d$, one can consider 
$\FRM\in W^{1,1}(I;L^1(\Omega;\R^d))$ and $\fRM\in W^{1,1}(I;L^1(\GN;\R^d))$
and then to perform the 
a-priori estimate \eqref{engr-ineq-disc} by using  the discrete
by-part integration (=summation) and the discrete Gronwall inequality,
and the coercivity \eqref{posit} instead of \eqref{posit-D}.
For the purpose of a-priori estimates only, one can also weaken 
\eqref{posit-D} to positive semi-definiteness (and in particular 
the rate-independent, inviscid problem with $\bbD=0$), 
although the convergence seems not guaranteed, cf.\ also 
Remark~\ref{rem-viscosity} below.
\end{remark}

\begin{remark}\label{rem-Dirichlet}
\upshape
If we assume time-dependent boundary conditions such that $u(t)=u_{\mbox{\tiny\rm D}}(t)$ a.e.~on $\GD$ for every $t\in [0,T]$ for some  $u_{\mbox{\tiny\rm D}}(t)\in W^{1,p}(\Omega;\R^d)$
Then the shift $u\mapsto u+u_{\mbox{\tiny\rm D}}(t)$ transform the problem to zero boundary conditions for $u$, i.e., $u_0=0$ on $\GD$.   

\end{remark}

\section{Computer implementation and illustrative simulations}\label{sec-simul}

We demonstrate varying mode-mixity of delamination on
a relatively simple two-dimensional example motivated by the pull-push shear 
experimental test used in engineering practice \cite{CoCa11ES}.
Intentionally, we use the same geometry, shown in Fig.~\ref{fig_m1},
as in \cite{RoMaPa13QMMD,RoPaMa??LSAQ} in order to have a comparison of our 
weak solution of the engineering non-associative visco-elastic model 
with a maximally-dissipative local solution and the energetic solution 
of the associative inviscid model presented respectively 
in \cite{RoPaMa??LSAQ} and in \cite{RoMaPa13QMMD}.
In contrast to Sections~\ref{sec2}--\ref{sec-ent-trans}, 
only one bulk domain is considered and $\GC$ is a part of its 
boundary but this modification is straightforward;
alternatively, one may also think about 
$\Omega_-$ as a completely rigid body in the previous setting.
Here $\Omega_+$ is a two-dimensional rectangular domain glued on the
most of its bottom side $\GC$ with 
the Dirichlet loading acting on the right-hand side $\GD$ in 
the direction $(1,0.6)$, cf.\  Fig.~\ref{fig_m1}, increasing linearly in time
with velocity 
$0.3\ $mm/s.

\begin{my-picture}{.95}{.17}{fig_m1}
\psfrag{GN}{\footnotesize $\GN$}
\psfrag{GD}{\footnotesize $\GD$}
\psfrag{GC}{\footnotesize $\GC$}
\psfrag{elastic}{\footnotesize elastic body}
\psfrag{obstacle}{\footnotesize rigid obstacle}
\psfrag{adhesive}{\footnotesize adhesive}
\psfrag{LC}{$L_c$}
\psfrag{L}{\footnotesize $L=\ $250\,mm}
\psfrag{H}{\footnotesize $H=$}
\psfrag{12.5}{\scriptsize 12.5\,mm}
\psfrag{loading}{\footnotesize loading}
\hspace*{-.5em}\vspace*{-.1em}{\includegraphics[width=.95\textwidth]{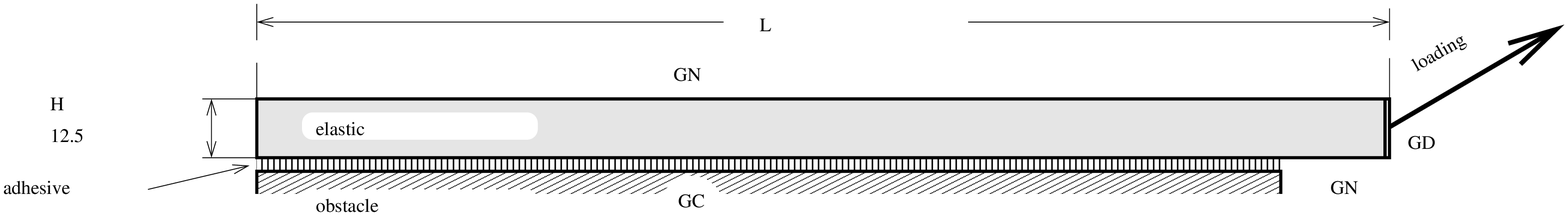}}
\end{my-picture}
\vspace*{-1.3em}
\begin{center}
{\small Fig.\,\ref{fig_m1}.\ }
\begin{minipage}[t]{.85\textwidth}\baselineskip=8pt
{\small
Geometry and boundary conditions of the problem considered.
The length of the initially glued part $\GC$ is $0.9L=225\,$mm, 
the adhesive layer has zero thickness.}
\end{minipage}
\end{center}

The bulk material is considered  linear, homogeneous, and isotropic
with the Young modulus $E=70$ GPa and Poisson's ratio $\nu=0.35$
(which corresponds to aluminum);
thus $\mathbb C_{ijkl}=\frac{\nu E}{(1{+}\nu)(1{-}2\nu)}\delta_{ij}\delta_{kl}
+\frac E{2(1{+}\nu)}(\delta_{ik}\delta_{jl}+\delta_{il}\delta_{jk})$
with $\delta_{ij}$ standing for the Kronecker symbol.
For the viscosity tensor we consider $\mathbb D=\chi\mathbb C$ with a relaxation time
$\chi=0.001\,$s, which is very small with the relation to the loading speed
we have considered; actually, we did not see any essential difference
for just merely elastic material $\chi=0$ although the model
has only a limited validity for this case, cf.\ Remark~\ref{rem-viscosity} below.

For the adhesive, we took a normal stiffness $\kappa_{\rm n}=$150 GPa/m, a tangential
stiffness with $\kappa_{\rm t}=\kappa_{\rm n}/2$, and the Mode-I fracture
toughness $a_{_{\rm I}}=187.5$ J/m$^2$. Furthermore,
the engineering model \eqref{a-mixity-dependence} was used 
with $\lambda=0.333$ (and with $\epsilon=0$)
which corresponds to a rather moderate 
mode-sensitivity $a(90^\circ)/a_{\rm I}=a_{\rm II}/a_{\rm I}=4$.

The numerical stability, i.e.\ the a-priori estimates \eqref{a-priori3}, 
fully applies.
The (rather negligible) shortcuts is that we take $p=2$ (instead of $p>2$)
so that the convergence in  Appendix below  applies
only up to this (small) discrepancy. 

For the computer implementation, it is important that \eqref{semi-impl+2} 
represents  a {\it linear-quadratic program},
which allows for algorithmically a very efficient solution. Furthermore,
\eqref{semi-impl+2} 
is even a {\it linear program}. On top of it, as we used P0-elements,
cf.\ \eqref{Zh}, it localizes on particular 
elements on $\GC$, so that its solution is trivial and fast;
for this effect, it is also important that we do not need to
consider any gradient of $z$, in contrast to the associative 
model \cite{tr+mk+jz,RoMaPa13QMMD,RoPaMa??LSAQ}.
For the Dirichlet loading, we have used Remark~\ref{rem-Dirichlet}.

A noteworthy attribute of our problem is that the inelastic process
of delamination occurs on the boundary $\GC$ while in the bulk domain(s) 
it is linear. This allows for elimination of nodal values arising 
by using P1-elements \eqref{Uh} inside $\Omega$ 
and considerable reduction of degrees of freedom by 
considering only nodal or element values on $\GC$.
In fact, this idea has been systematically exploited even on the 
continuous level when implementing the boundary-element method, cf.\ 
\cite{Mukh01BIEC,ParCan97BEMF,RoMaPa13QMMD,RoPaMA13QACV,RoPaMa??LSAQ,TMGP11BEMA},
although it is still not fully supported by a convergence analysis like 
Proposition~\ref{th:3.0} below due to general substantial 
theoretical difficulties related to this method.

Anyhow, for the computational experiments presented here partially also
with the goal to document the usually not investigated modelling issues
(in particular the energetics), we use a shortcut in 
implementing the spacial discretization \eqref{semi-impl+} with \eqref{Uh-Vh}
and, instead of an algebraic elimination of the interior nodal
point, we made this elimination by the collocation boundary-element method. 
Here also our choice $\mathbb D=\chi\mathbb C$ makes the implementation
of this method simpler, cf.~\cite[Remark 6.2]{Roub13ACVE}.
For the results presented on Figures~\ref{fig_m2}-\ref{fig_m3}, we have
used 81 elements on $\GC$, i.e.\ $h=2.7\bar{7}\,$mm
(=the size of a boundary segment in equidistant discretization), and 
the time step $\tau = 2.2\bar{2}\,$ms. 

This example exhibits remarkably varying mode of delamination. At the beginning 
the delamination is performed by a mixed mode close to Mode I given
essentially by the direction of the Dirichlet loading, cf.\ Figure~\ref{fig_m1},
while later it turns rather to nearly pure Mode II. Yet, at the 
very end of the process, due to elastic bending
the delamination starts performing also from the 
left-hand side of the bar opposite to the loading side, and
thus again a mixed mode occurs. This relatively complicated mixed-mode behaviour 
is depicted in Figure~\ref{fig_m2}(right), showing essential qualitative 
difference from the energetic solution which exhibits a non-physical 
tendency to slide to less-dissipative Mode I, cf.\ \cite[Fig.\,7]{RoMaPa13QMMD}.
We have here the energy 
(im)balance \eqref{total-energy} as an important ingredient and thus,
in contrast to usual engineering calculations, we trace also the 
energetics of the model, depicted on Figure~\ref{fig_m2}(left).

\begin{my-picture}{.95}{.4}{fig_m2}
\hspace*{-1.9em}\vspace*{-.1em}{\includegraphics[width=.59\textwidth,height=.4\textwidth]{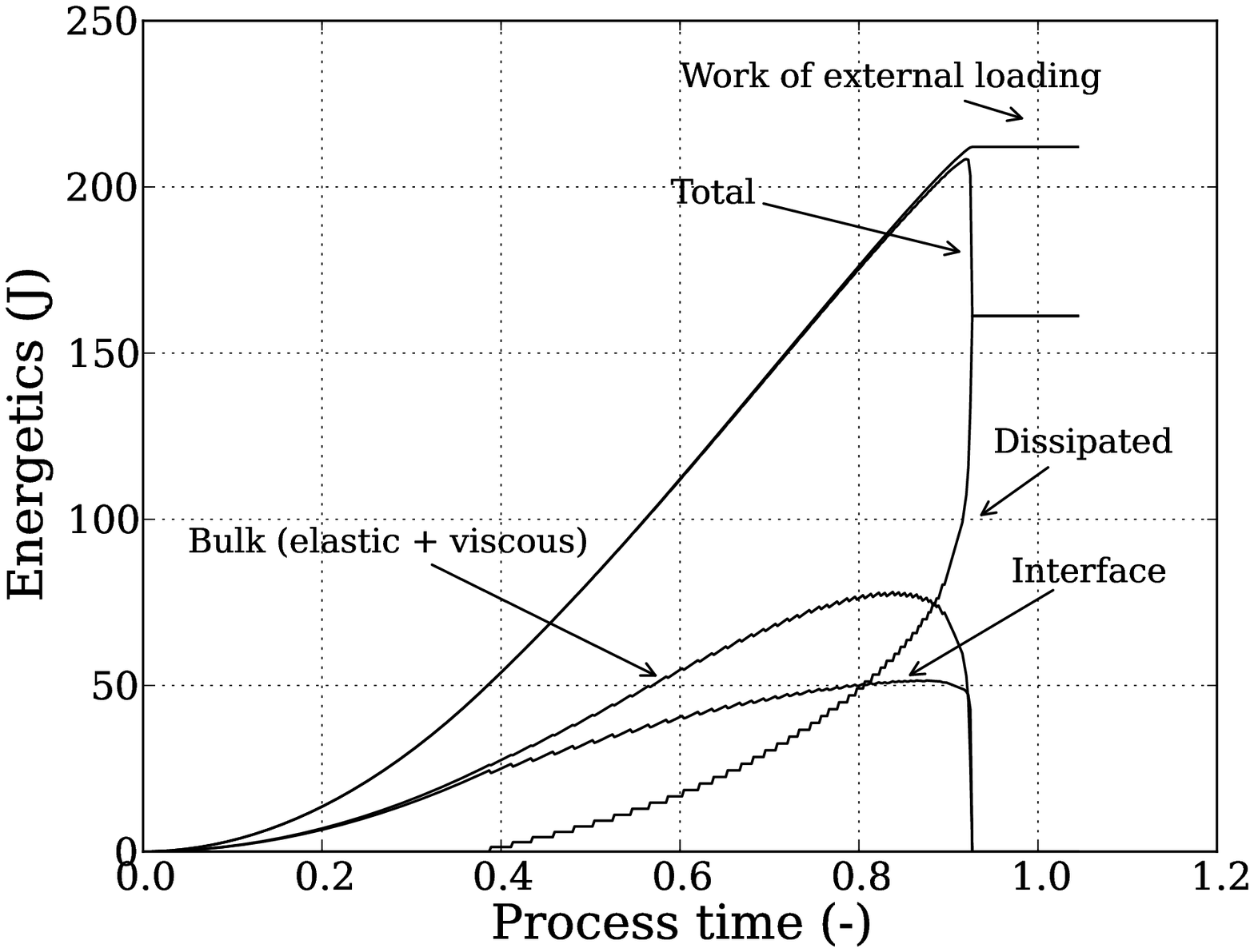}}
\hspace*{-2.em}\vspace*{-.1em}{\includegraphics[width=.49\textwidth,height=.4\textwidth]{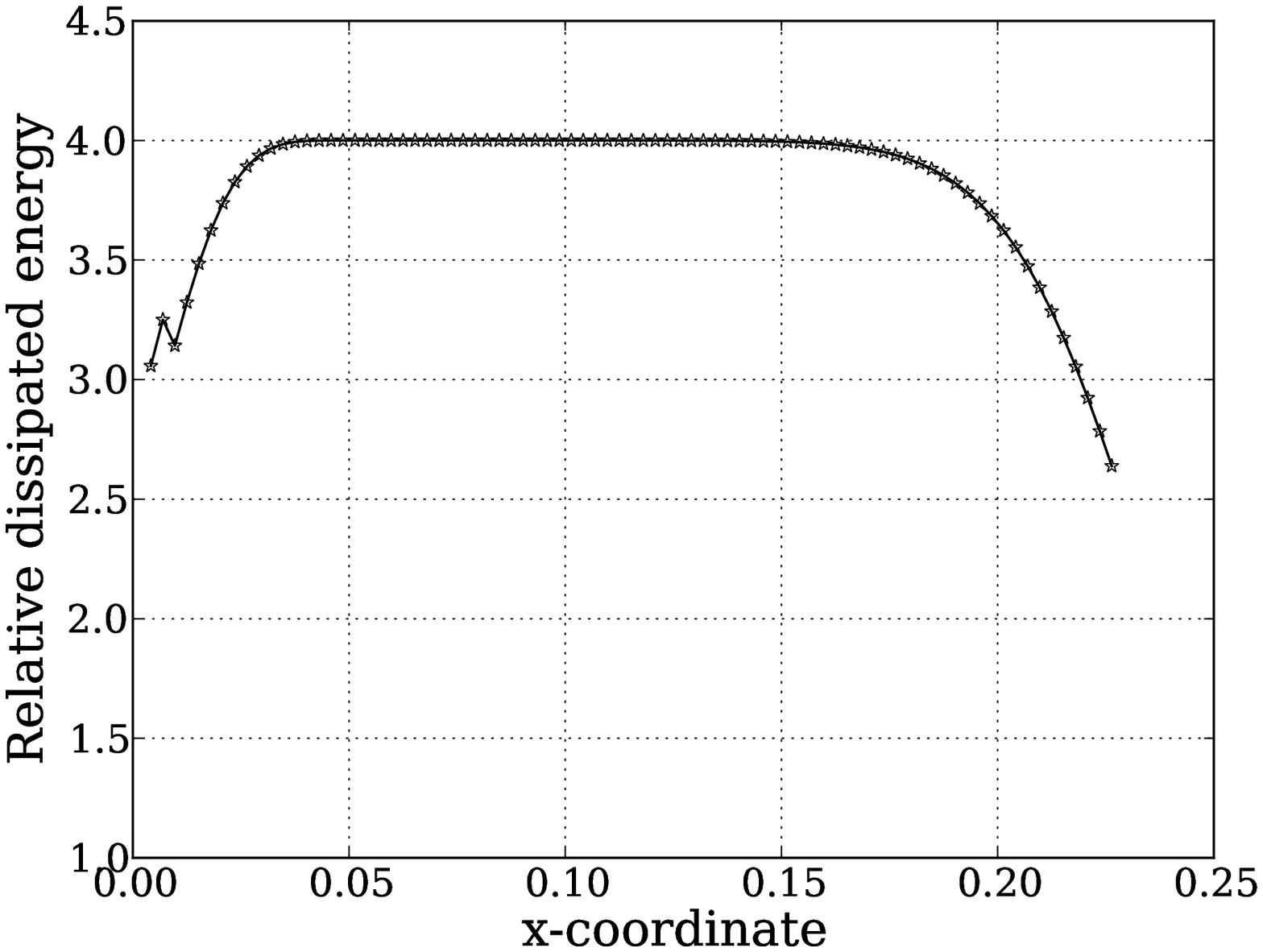}}
\end{my-picture}
\vspace*{-1em}
\begin{center}
{\small Fig.\,\ref{fig_m2}.\ }
\begin{minipage}[t]{.88\textwidth}\baselineskip=8pt
{\small
{\sf Left:} Time evolution of the energies:
the bulk and the interfacial parts of 
$\calE(t,\bar u_\tau(t),\bar z_\tau(t))+\int_0^t\calR_2(\DT u_\tau)\d t$,
the interfacial dissipated energy $\calR_1(\bar z_\tau(t){-}z_0)$,
their sum=total energy (=the left-hand side of \eqref{total-energy}),
and the work of external loading 
(=the right-hand side of \eqref{total-energy}).
\\
{\sf Right:} Mixity-mode distribution along $\GC$
evaluated according the overall dissipated energy related to 
$a_{_{\rm I}}$ after the delamination has been completed: 
value\,=\,1$\,\sim\,$Mode I, value\,=\,4\,=$\,a_{_{\rm II}}/a_{_{\rm I}}\!\sim\,$Mode II.

}
\end{minipage}
\end{center}

The evolution of the deformation of the visco-elastic domain (through the displacement $u$) 
and spacial distribution of the delamination parameter $z$ are depicted in 
Figure~\ref{fig_m3} at eight snapshots selected not uniformly to
visualize interesting effects when delamination starts to be completed. 
In particular, the delamination propagating from both sides at the very end 
(mentioned already above) is seen there.



\vspace*{-17em}
\begin{my-picture}{.95}{.95}{fig_m3}
\begin{tabular}{lll}
\hspace*{.3em}\LARGE $^{^{^{^{\mbox{\footnotesize $k$=100}}}}}$  & \hspace*{-.5em}\vspace*{-.1em}$^{^{^{^{^{^{{\includegraphics[width=.6\textwidth]{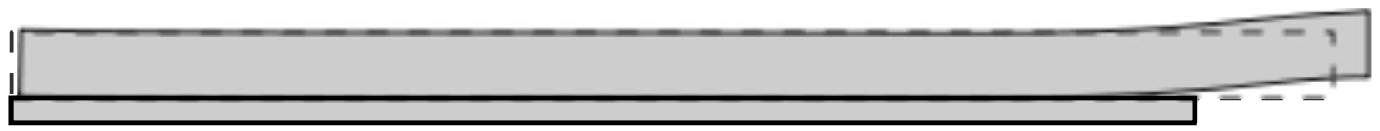}}}}}}}}$&
\hspace*{-.5em}{\includegraphics[width=.2\textwidth,height=.12\textwidth]{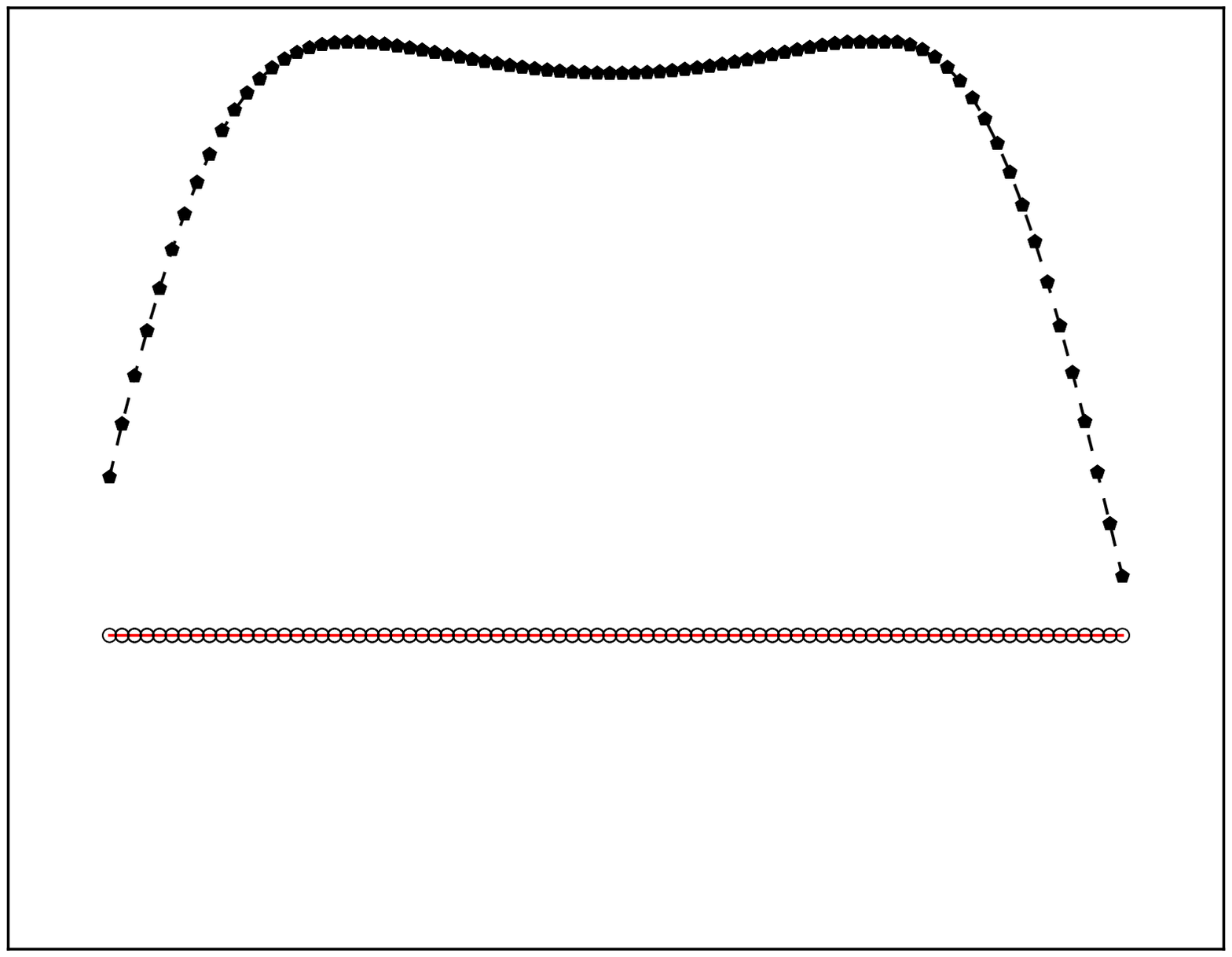}}
\\
\hspace*{.3em}\LARGE $^{^{^{^{\mbox{\footnotesize $k$=180}}}}}$  & \hspace*{-.5em}\vspace*{-.1em}$^{^{^{^{^{^{{\includegraphics[width=.6\textwidth]{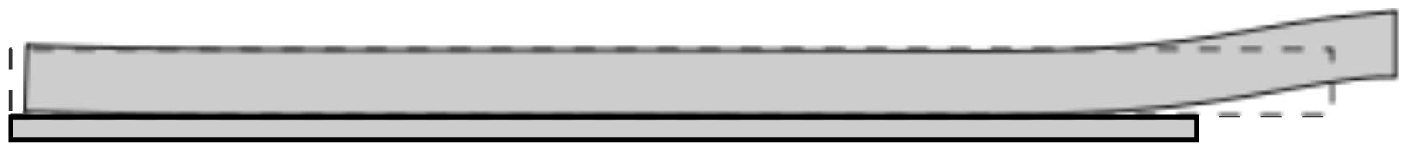}}}}}}}}$&
\hspace*{-.5em}{\includegraphics[width=.2\textwidth,height=.12\textwidth]{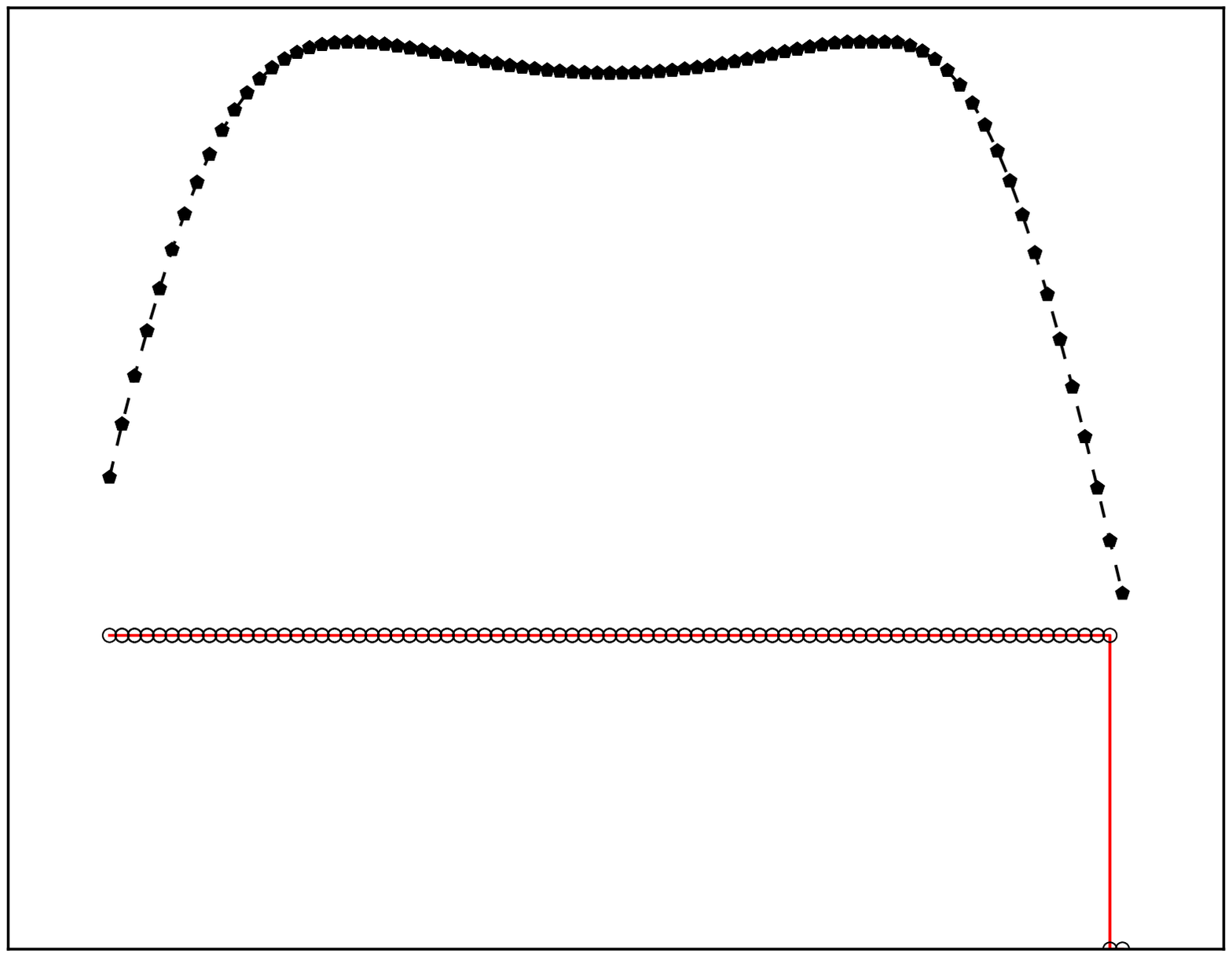}}
\\
\hspace*{.3em}\LARGE $^{^{^{^{\mbox{\footnotesize $k$=260}}}}}$  & \hspace*{-.5em}\vspace*{-.1em}$^{^{^{^{^{^{{\includegraphics[width=.6\textwidth]{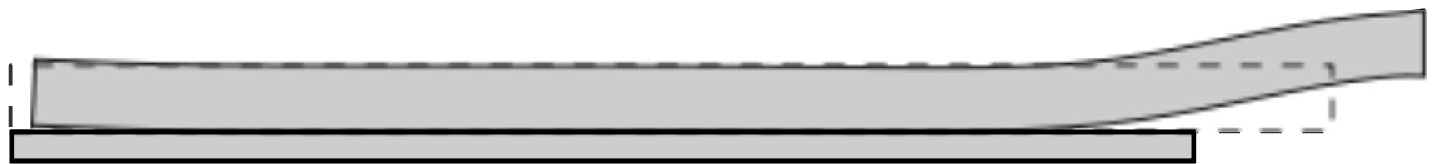}}}}}}}}$&
\hspace*{-.5em}{\includegraphics[width=.2\textwidth,height=.12\textwidth]{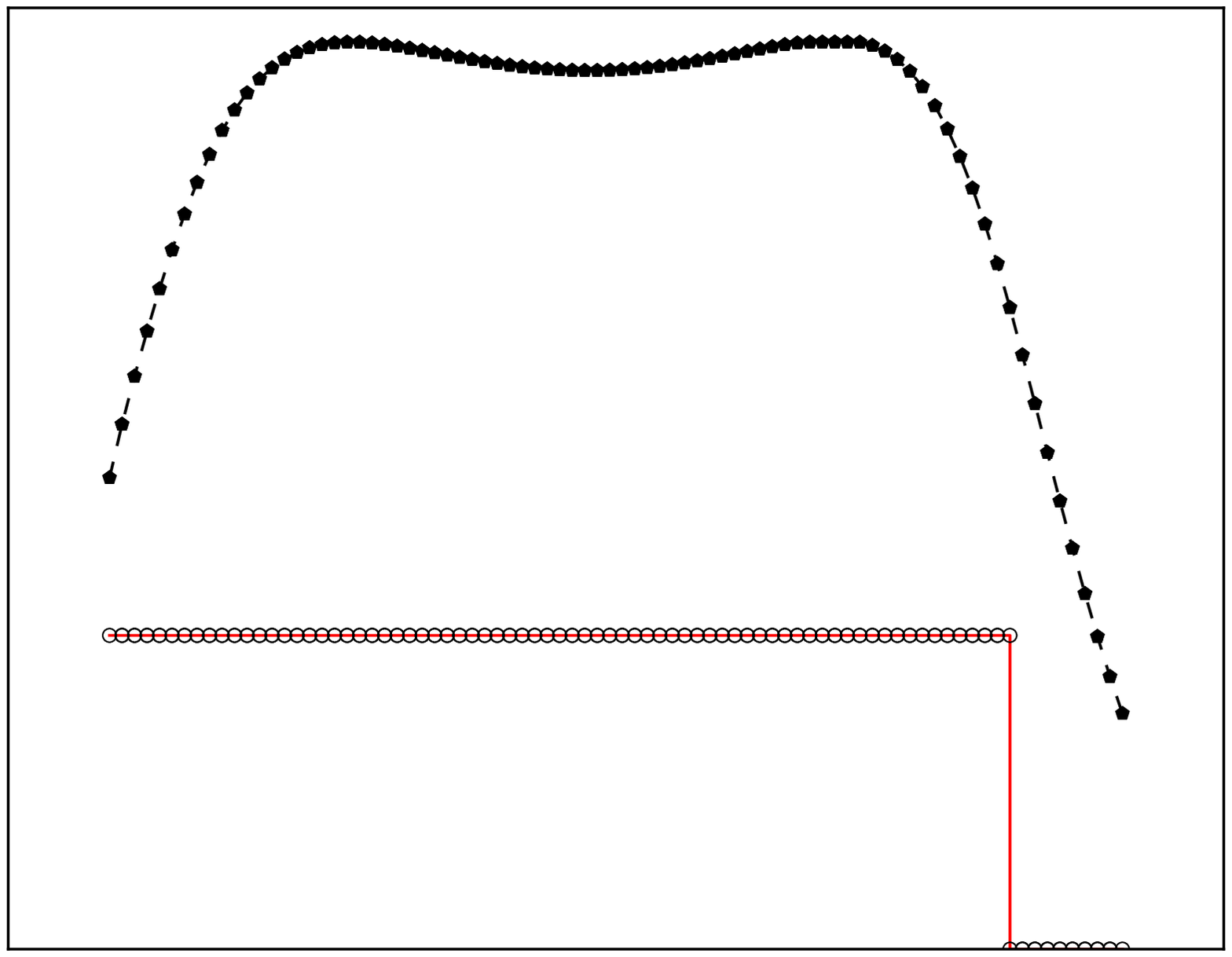}}
\\
\hspace*{.3em}\LARGE $^{^{^{^{\mbox{\footnotesize $k$=330}}}}}$  & \hspace*{-.5em}\vspace*{-.1em}$^{^{^{^{^{^{{\includegraphics[width=.6\textwidth]{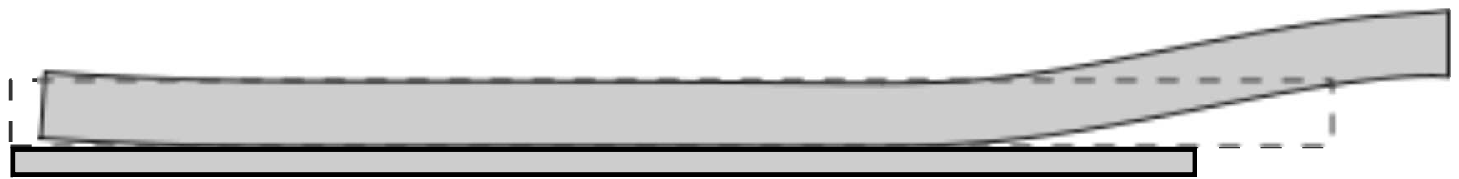}}}}}}}}$&
\hspace*{-.5em}{\includegraphics[width=.2\textwidth,height=.12\textwidth]{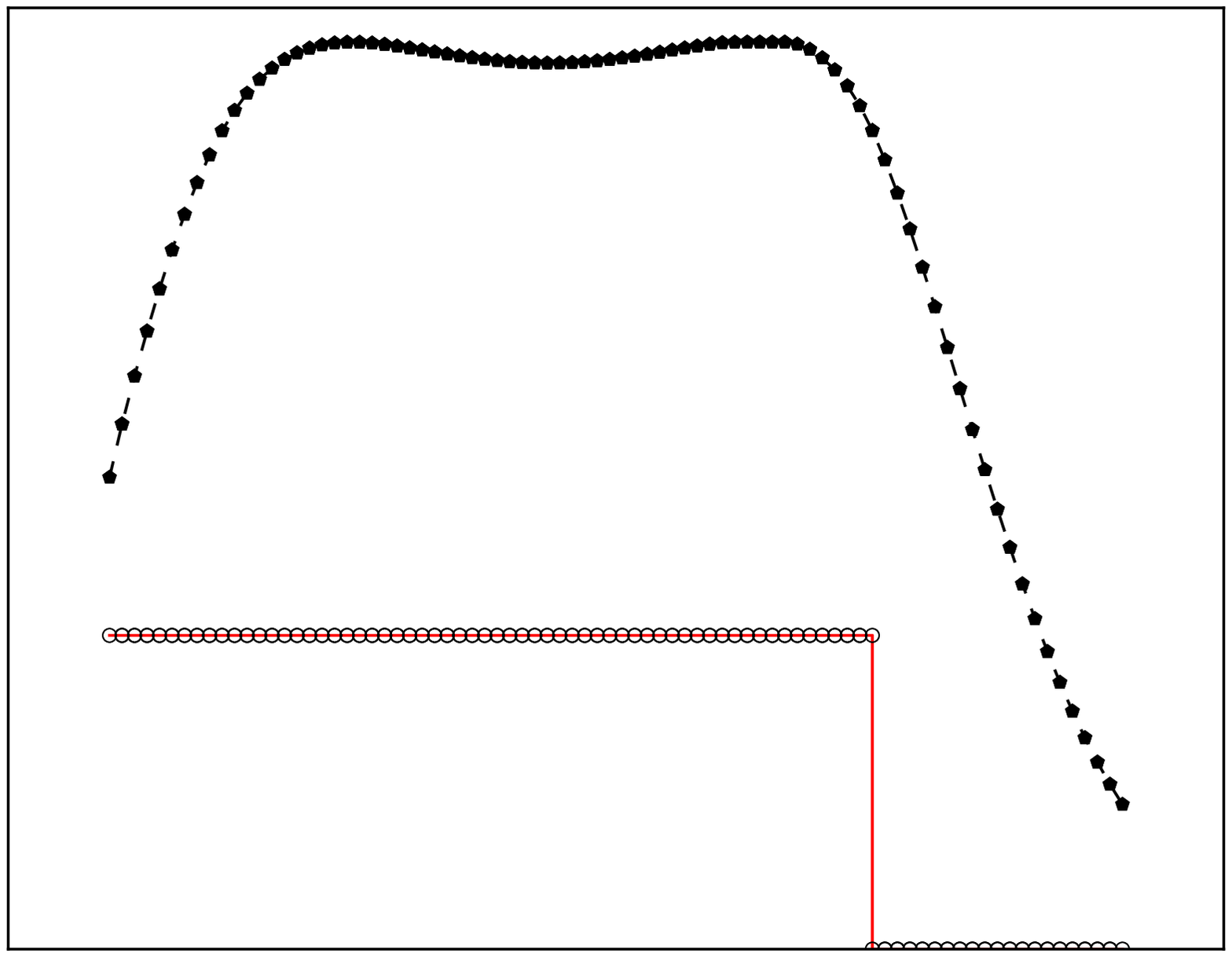}}
\\
\hspace*{.3em}\LARGE $^{^{^{^{\mbox{\footnotesize $k$=390}}}}}$  & \hspace*{-.5em}\vspace*{-.1em}$^{^{^{^{{\includegraphics[width=.6\textwidth]{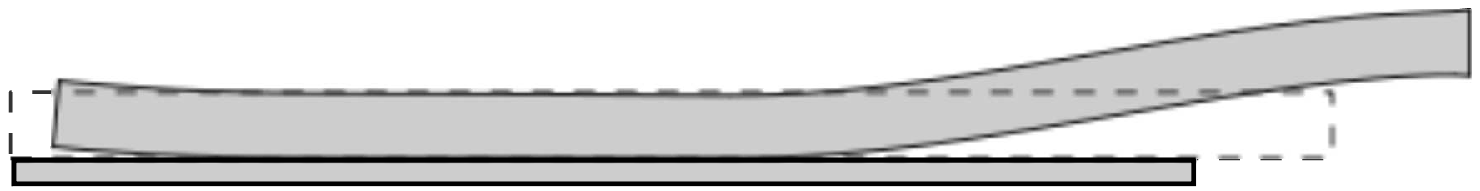}}}}}}$&
\hspace*{-.5em}{\includegraphics[width=.2\textwidth,height=.12\textwidth]{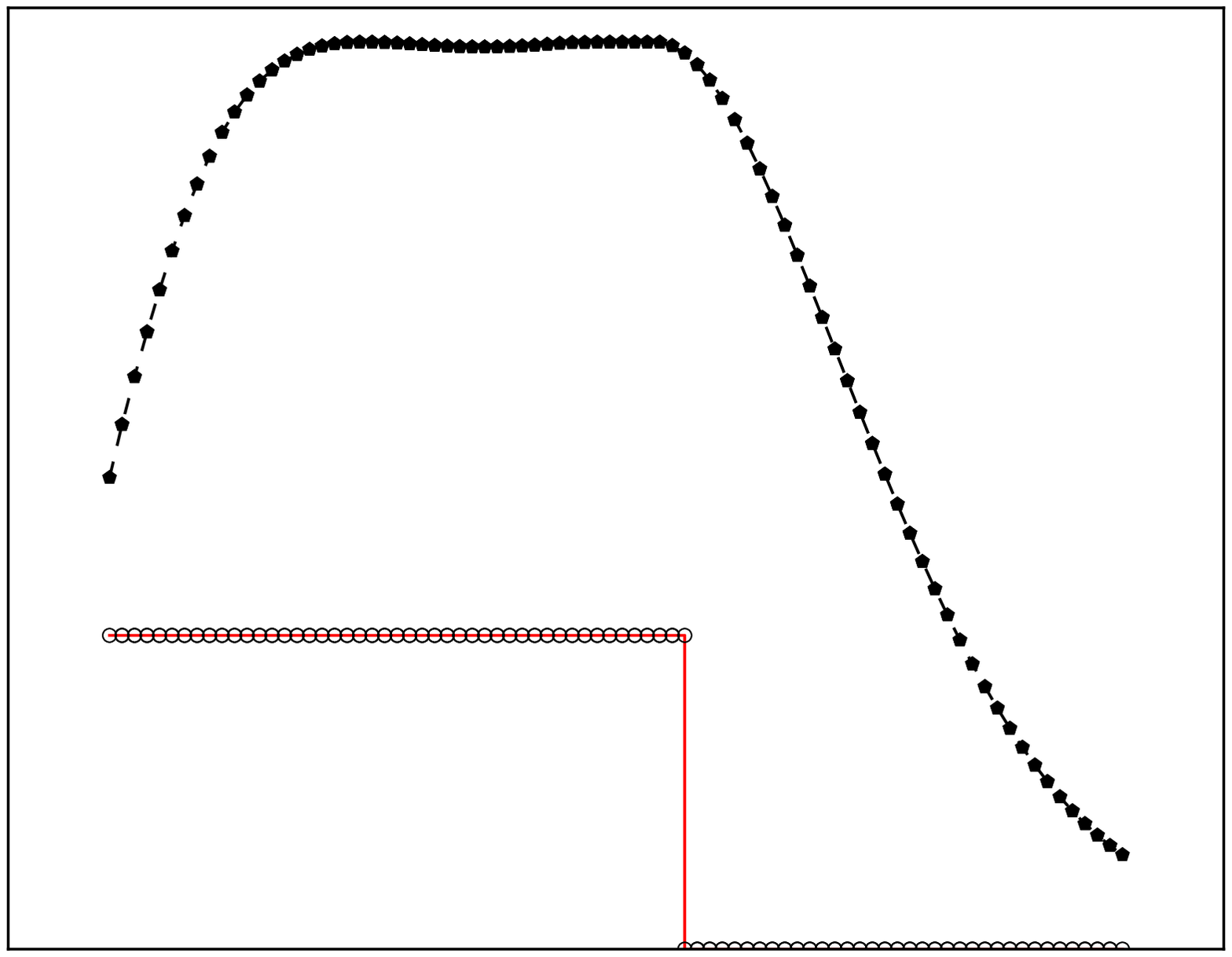}}
\\
\hspace*{.3em}\LARGE $^{^{^{^{\mbox{\footnotesize $k$=412}}}}}$  & \hspace*{-.5em}\vspace*{-.1em}$^{^{^{^{{\includegraphics[width=.6\textwidth]{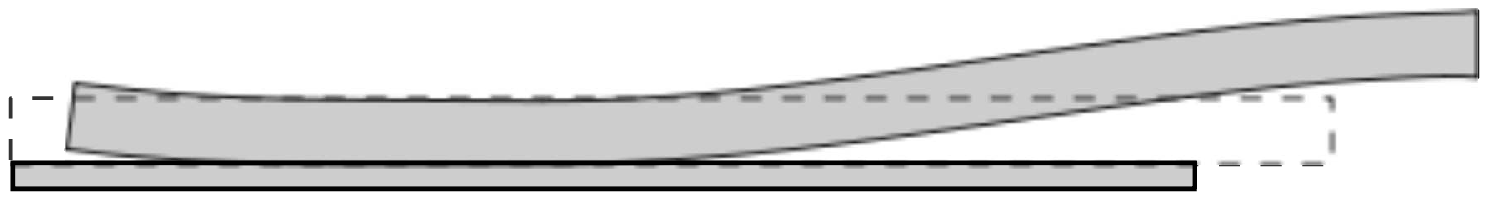}}}}}}$&
\hspace*{-.5em}{\includegraphics[width=.2\textwidth,height=.12\textwidth]{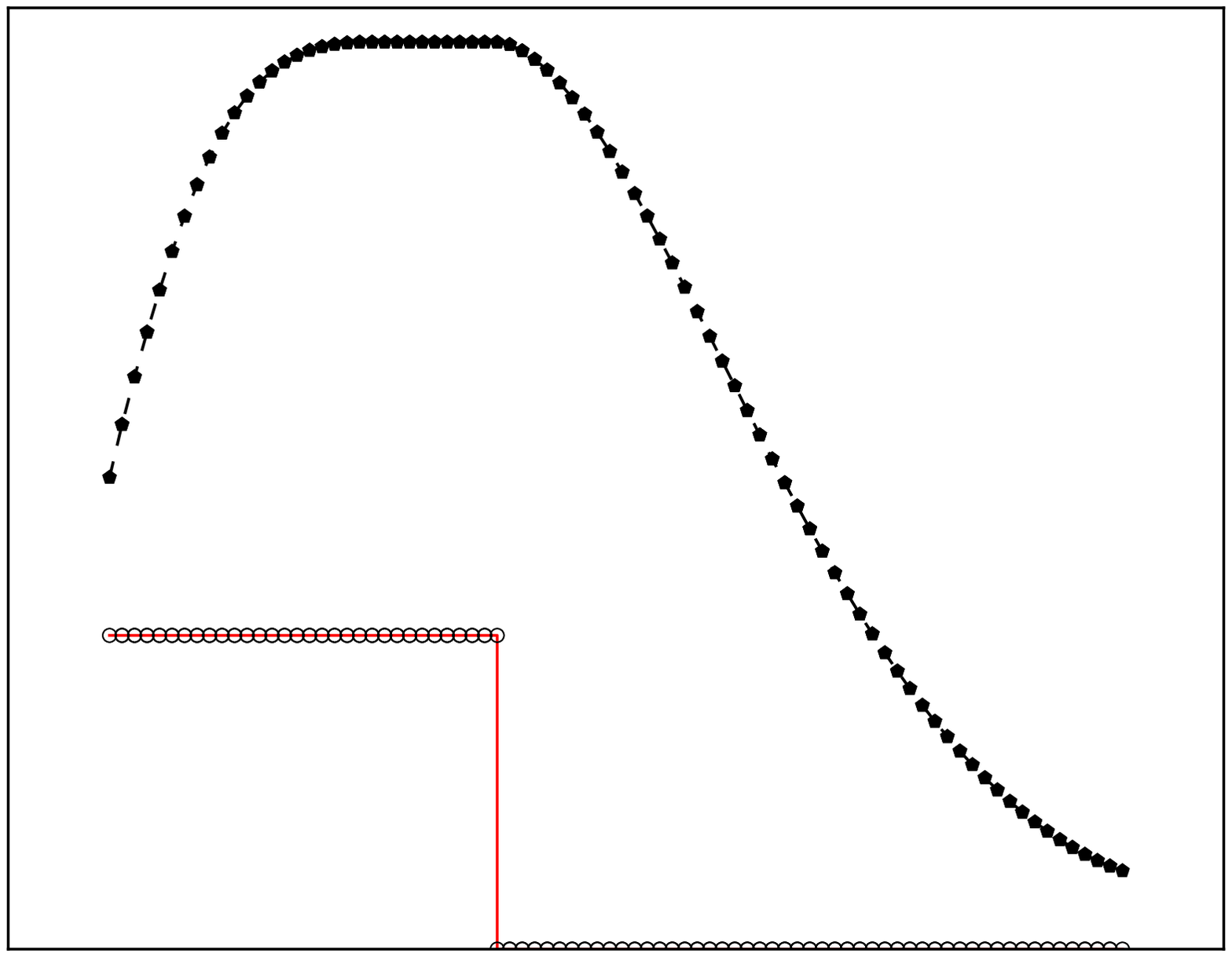}}
\\
\hspace*{.3em}\LARGE $^{^{^{^{\mbox{\footnotesize $k$=416}}}}}$  & \hspace*{-.5em}\vspace*{-.1em}$^{^{^{^{{\includegraphics[width=.6\textwidth]{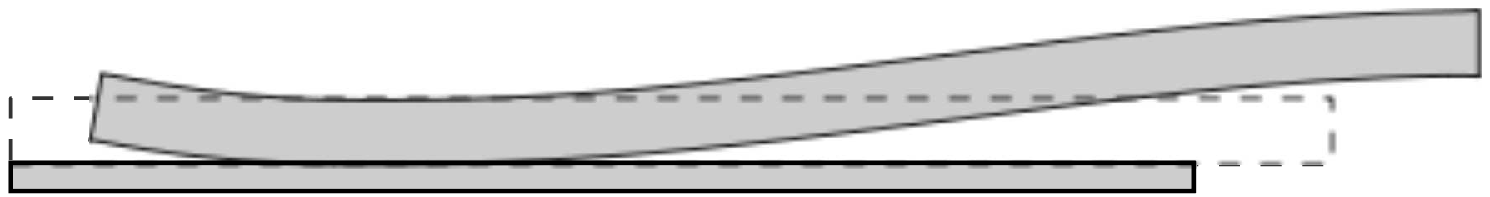}}}}}}$&
\hspace*{-.5em}{\includegraphics[width=.2\textwidth,height=.12\textwidth]{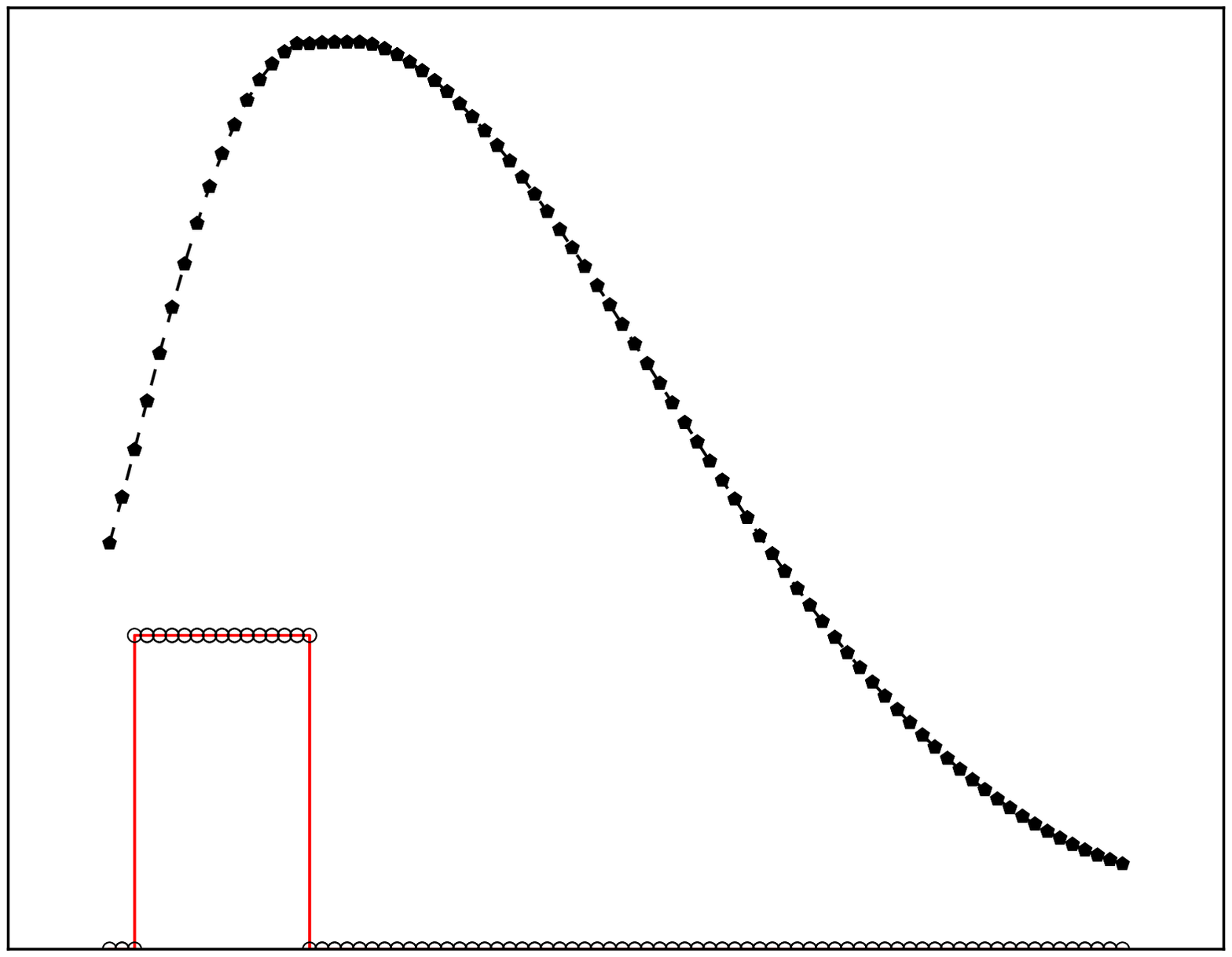}}
\\
\hspace*{.3em}\LARGE $^{^{^{^{\mbox{\footnotesize $k$=417}}}}}$  & \hspace*{-.5em}\vspace*{-.1em}$^{^{{\includegraphics[width=.6\textwidth]{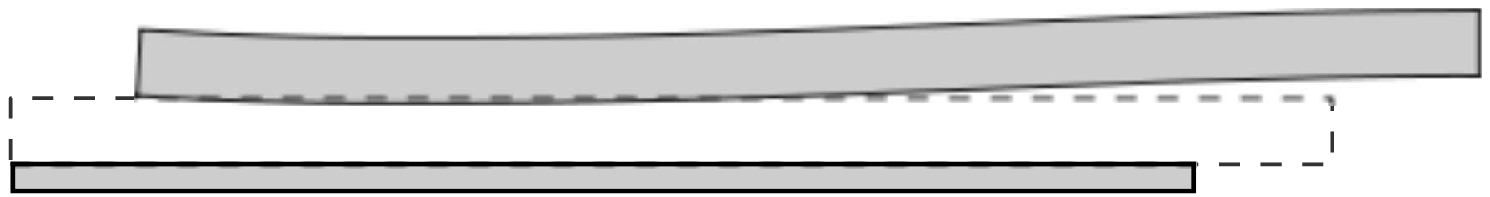}}}}$&
\hspace*{-.5em}{\includegraphics[width=.2\textwidth,height=.12\textwidth]{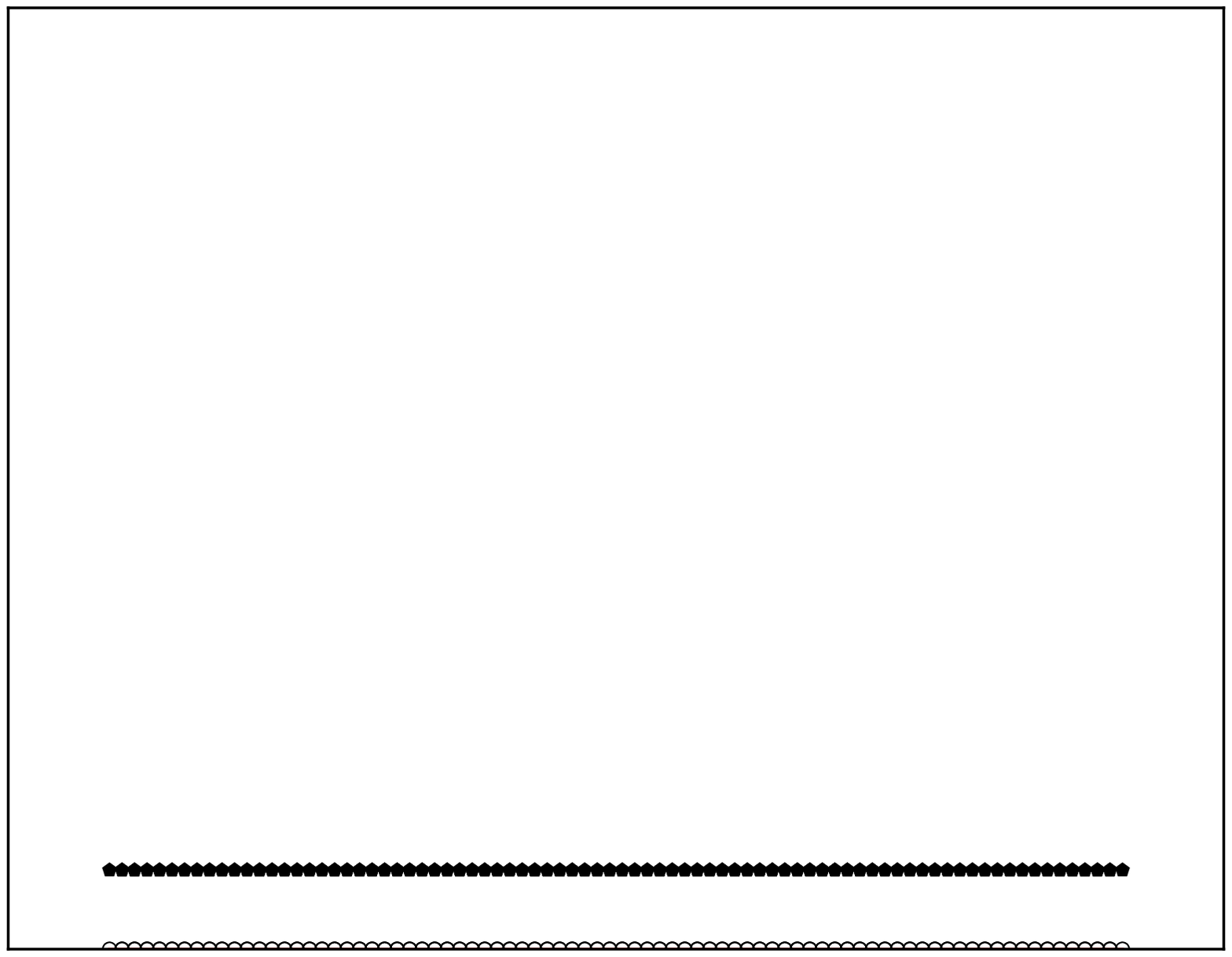}}
\end{tabular}
\end{my-picture}
\vspace{20em}
\begin{center}
{\small Fig.\,\ref{fig_m3}.\,}
\begin{minipage}[t]{.89\textwidth}\baselineskip=8pt
{\small
Time evolution at eight snapshots of the geometrical configuration
until complete delamination (displacement depicted magnified $100\,\times$)
and the spacial distribution of the mode-mixity angle $\psiG$ and the delamination
parameter $z$ along the interface $\GC$.
}
\end{minipage}
\end{center}

\medskip

Eventually, the joint convergence from Proposition~\ref{th:3.0} below
for time- and FEM-spatial discretization
(although here implemented by BEM) is demonstrated in Figures~\ref{fig_conv}
and \ref{fig_conv+} for two different gradually refining 
space/time discretizations. 
We choose the scenario keeping the ratio $\tau/h$ constant, although 
Proposition~\ref{th:3.0} itself does not give any particular suggestion
in this respect. Anyhow, the tendency of 
convergences is clearly seen, although we naturally do not know the 
exact solution so that we cannot evaluate any actual error. On top of it,
the exact solution does not need to be unique so we even do not have guaranteed 
the convergence of the whole sequence of the approximate solutions and, moreover,
the simplified implementation by 
collocation BEM does not have guaranteed convergence,
in contrast to FEM stated in Proposition~\ref{th:3.0}.



\begin{my-picture}{.95}{.38}{fig_conv}
\hspace*{-1.em}\vspace*{-.1em}{\includegraphics[width=.5\textwidth]{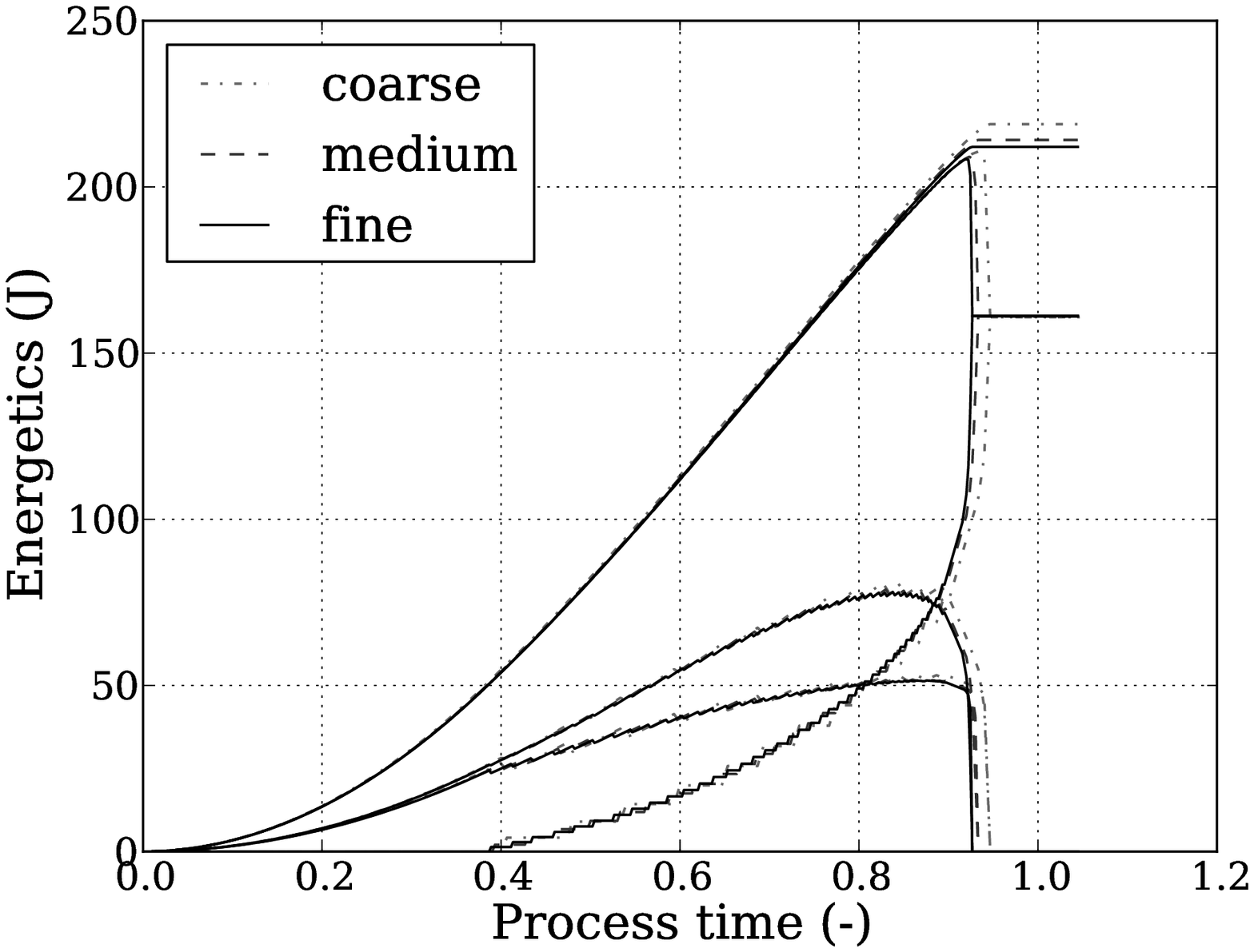}}
\hspace*{-.5em}\vspace*{-.1em}{\includegraphics[width=.5\textwidth]{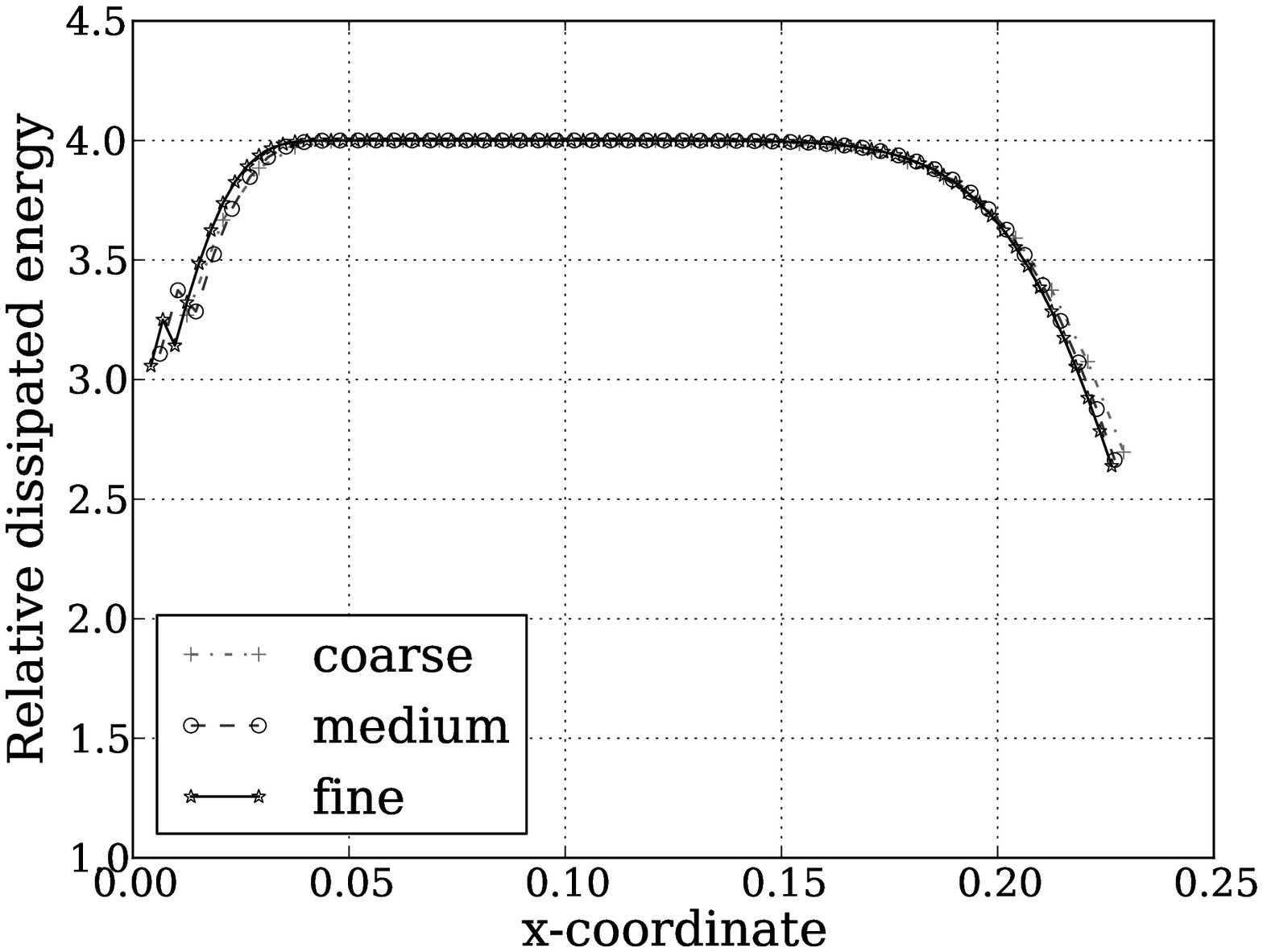}}
\end{my-picture}
\vspace{-1em}
\begin{center}
{\small Fig.\,\ref{fig_conv}.\ }
\begin{minipage}[t]{.8\textwidth}\baselineskip=8pt
{\small
Convergence test of the quantities from 
Figure~\ref{fig_m2}.  
Coarse\,=\,27 elements on $\GC$, 
medium\,=\,54 elements, and fine\,=\,81 elements. The ratio $\tau/h$ is constant.
}
\end{minipage}
\end{center}

\begin{my-picture}{.95}{.35}{fig_conv+}
\hspace*{-1.em}\vspace*{-.1em}{\includegraphics[width=.5\textwidth]{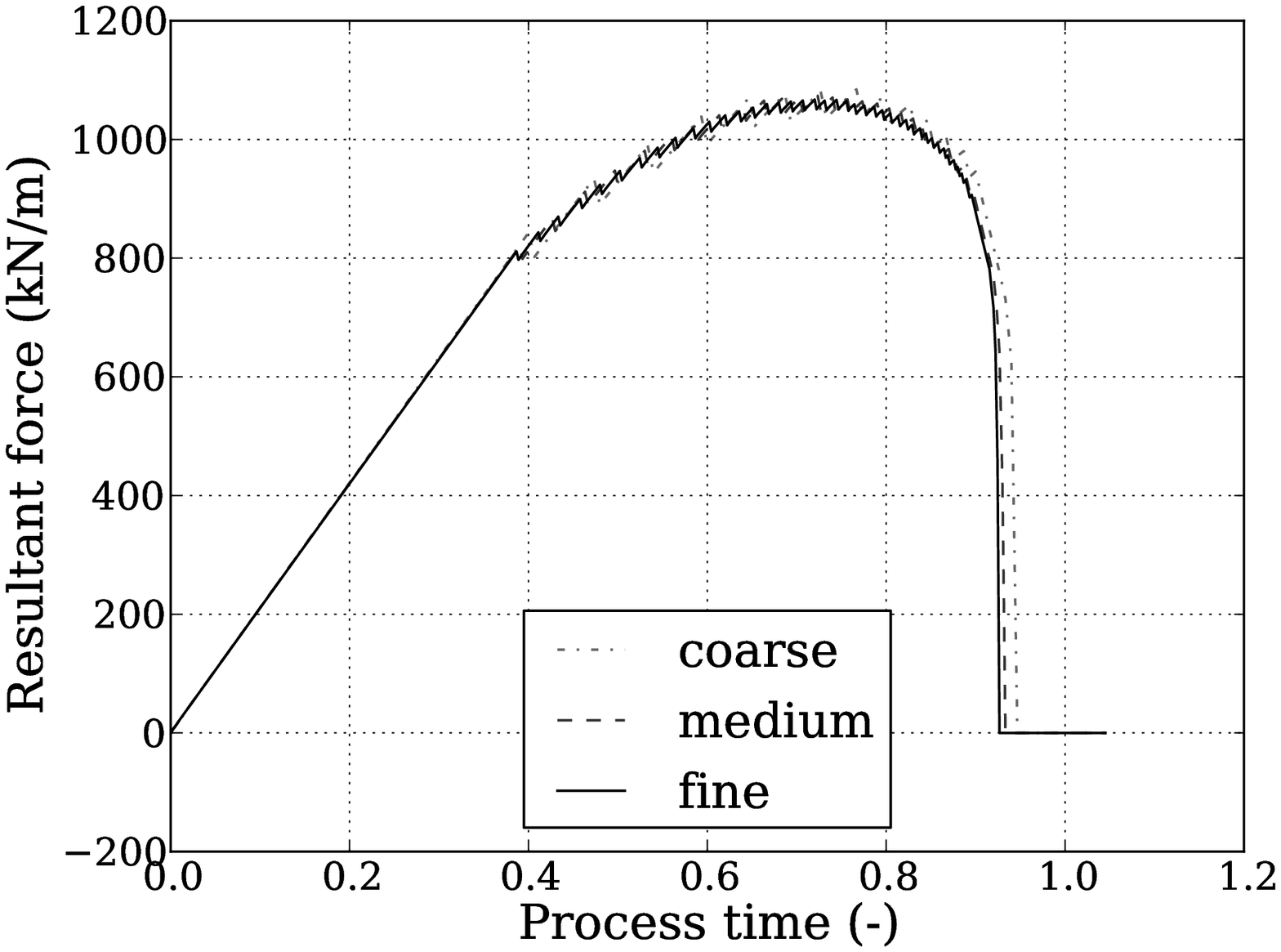}}
\hspace*{-.5em}\vspace*{-.1em}{\includegraphics[width=.5\textwidth]{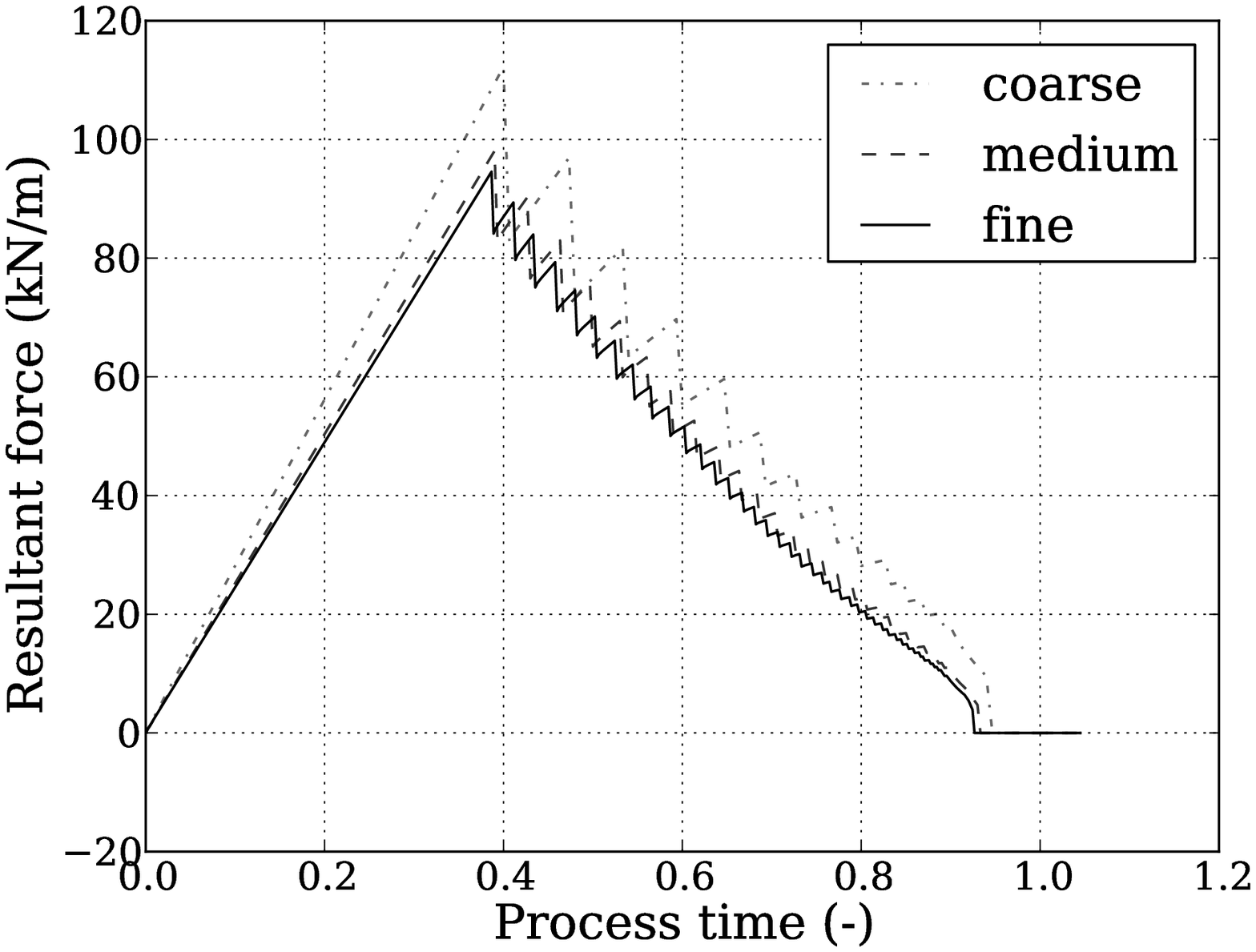}}
\end{my-picture}
\vspace{-1em}
\begin{center}
{\small Fig.\,\ref{fig_conv+}.\ }
\begin{minipage}[t]{.8\textwidth}\baselineskip=8pt
{\small
Convergence test of the total force response evolving in time: the horizontal 
(left) and the vertical (right) components.
}
\end{minipage}
\end{center}

\begin{remark}[{\it Non-conservation of energy}]\label{rem-engr}
\upshape
It is interesting to check the energy (im)balance \eqref{total-energy}.
In Figure\,\ref{fig_m2}(left), we can see it depicted for $t_1=0$ as a function of
time $t_2$: the upper line is the right-hand side of  \eqref{total-energy}
while the line below is the left-hand side of \eqref{total-energy}.
We can clearly see that the difference is not zero and is increasing in time, which
is in accord with \eqref{total-energy} because otherwise,
if the difference would decrease on some time interval $[t_1,t_2]$, 
\eqref{total-energy} could not be valid on this interval.
This non-vanishing difference between the left- and the right-hand sides of 
\eqref{total-energy} is likely not because of a possible 
 numerical possible error (as Fig.\,\ref{fig_conv} shows a nice convergence)
but seems to have a physical meaning that some part of energy is lost 
(dissipated) due to some neglected dissipation mechanisms which
would guarantee energy conservation. This would certainly happen if 
the bulk viscosity would be completely neglected, as now well 
documented in \cite{Roub13ACVE,RoPaMA13QACV}. Here, although 
we have bulk viscosity, this viscosity is 
of a slower growth  
than the stored energy and thus might not
be strong enough for fast processes, as observed at the fast completion
of the delamination process here. Nevertheless, superquadratic growth of the viscosity 
would recast our model to the so-called doubly-nonlinear problem and would bring additional mathematical difficulties.  
\end{remark}

\begin{remark}[{\it Comparison with other models/solutions}]\label{rem-cmp}
\upshape
As the viscosity has been chosen low, this example bears a good qualitative 
comparison (including the energy gap discussed in Remark~\ref{rem-engr})
with the associative inviscid model with interfacial plasticity 
from \cite{tr+mk+jz,RoMaPa13QMMD}, provided maximally-dissipative local 
solutions are considered \cite{RoPaMa??LSAQ}.
The non-associative model allows for a general monotone dependence $a(\cdot)$
not necessarily just the special ansatz \eqref{a-mixity-dependence}
which enables a possible fitting to that associative model 
which involves only few parameters. Also, 
the non-associative model allows do not involve any gradients of $z$
which allows simpler computational implementation, as mentioned already 
above. On the other hand, the ``vanishing-viscosity'' asymptotics of 
the non-associative model towards the inviscid, rate-independent model
is not entirely clear, cf.~Remark~\ref{rem-viscosity} below.
\end{remark}
\section{Appendix: convergence}\label{sec-main-res}


In terms of the piecewise constant/affine interpolants introduced in 
Section~\ref{sec-ent-trans}, cf.\ \eqref{interpolants}, the discrete scheme 
\eqref{disc} summed for $k=1,...,T/\tau$ can be written ``more compactly'' as:
\begin{subequations}\label{disc+}
\begin{align}\nonumber
 &\forall v\!\in\!L^1(I;V_h):\ \ 
\int_{Q\setminus\SC}\!\!\Big(\bbD e\big(\DT u_{\tau h}\big)
{+}\bbC(e(\bar u_{\tau h}))\Big){:}e(\testu{-}\bar u_{\tau h})\,\d x\d t
\\[-.5em]\nonumber&\hspace{10em}
+\int_{\SC}\!\!\!\underline z_{\tau h}\bbA\JUMP{\bar u_{\tau h}}
{\cdot}\JUMP{\testu{-}\bar u_{\tau h}}\d S\d t
\\[-.5em]\label{disc-u+}
&\hspace{15em}
\ge\int_Q\!\bar \FRM_\tau{\cdot}(\testu{-}\bar u_{\tau h})\,\d x\d t
+\int_{\Sigma_N}\!\!\bar \fRM_\tau{\cdot}(\testu{-}\bar u_{\tau h})\,\d S\d t
\\&\nonumber
\forall t_1\!<\!t_2\,,\ t_1,t_2\!\in\!\{k\tau\}_{k=1}^{T/\tau}:\ \ \ \
\Phi\big(u_{\tau h}(t_2),z_{\tau h}(t_2)\big)
+\int_{t_1}^{t_2}\!
\calR_1\big(\bar u_{\tau h};\DT z_{\tau h}\big)+\calR_2\big(\DT u_{\tau h}\big)
\,\d t
\\[-.5em]\label{disc-z-1+}&\hspace{8em}
\le\Phi\big(u_{\tau h}(t_1),z_{\tau h}(t_1)\big)
+\int_{t_1}^{t_2}\!\!\int_\Omega\bar\FRM_\tau{\cdot}\DT u_{\tau h}\,\d x\d t
+\int_{t_1}^{t_2}\!\!\int_{\GN}\!\!\bar\fRM_\tau{\cdot}\DT u_{\tau h}\,\d S\d t,
\\&\label{disc-z-2+}
\forall \tilde{z}\!\in\!Z_h:\qquad \Phi\big(\bar u_{\tau h}(t),\bar z_{\tau h}(t)\big)
\le\Phi\big(\bar u_{\tau h}(t),\tilde z\big) 
+\calD\big(\bar u_{\tau h}(t);0,\tilde z-\bar z_{\tau h}(t)\big).
\end{align}
\end{subequations}

We now enlist further assumptions needed for convergence:
\begin{subequations}\label{ass-for-conv}\begin{align}\label{ass-monot-phi}
&\exists\eta>0\ \forall e_1,e_2\in\R_{\mathrm{sym}}^{d\times d}:\qquad
(\varphi'(e_1)-\varphi'(e_2)){:}(e_1-e_2)\ge
\eta|e_1{-}e_2|^p,
\\\label{anot1}
&\alpha \in C (\R^d;\R), 
\\\label{ass-mesh}
&\bigcup_{h>0} V_h\text{ dense in }W^{1,p}(\Omega{\setminus}\GC;\R^d)\ \&\  V_{h_0}\subset V_{h_1} \mbox{ if } h_1<h_0\ .
\end{align}\end{subequations}


\begin{proposition}[Unconditional convergence]\label{th:3.0}
Let the assumptions of Proposition~\ref{prop:apriori} with $p>d$ as well as 
\eqref{ass-for-conv} be fulfilled. Then, for $\tau\to0$ and $h\to0$, the 
approximate solutions $(u_{\tau h},z_{\tau h})$ converge in terms of subsequences
in the sense that
\begin{subequations}\begin{align}\label{conv-u}
&&&u_{\tau h}\to u&&\text{in }L^q(I;W^{1,p}(\Omega{\setminus}\GC;\R^d)),
&&&&
\\\label{conv-z}
&&&\bar z_{\tau h}(t)\weaksto z(t)
&&\text{ in }L^\infty(\GC)\qquad\text{for all }t\in\bar I&&&&
\end{align}
\end{subequations}
for any $1\le q<\infty$, and any $(u,z)$ obtained by such converging 
subsequences is a weak solution to the adhesive contact problem 
due to Definition~\ref{def4}.
In particular, such weak solutions do exist.
\end{proposition}
%

\noindent{\it Sketch of the proof.}
By the a-priori estimates \eqref{a-priori3} and by Banach's
principle, we can select a weakly* converging subsequence in the
spaces indicated in \eqref{a-priori3}. Moreover, by Helly's selection
principle, we can consider this subsequence so that \eqref{conv-z} holds.

Now we improve  the weak* convergence $u_{\tau h}\weaksto u$ 
in the space $L^\infty(I;W^{1,p}(\Omega{\setminus}\GC;\R^d))\cap H^1(I;H^1(\Omega{\setminus}\GC;\R^d))$
by proving also  
the strong convergence in $L^p(I;W^{1,p}(\Omega{\setminus}\GC;\R^d))$. 
To this goal, let us 
exploit the uniform monotonicity of the operator $-{\rm div}\,
\bbC(e(\cdot))$, 
cf.\ \eqref{ass-monot-phi}, and use \eqref{disc-u+} for 
$v=2\bar u_{\tau h}{-}v_h$ with some $v_h\in L^\infty(I;V_h)$:
\begin{align}\nonumber
&\!\!\!\eta\|e(\bar u_{\tau h}{-}v_h)\|_{L^p(Q;\R^{d\times d})}^p
\le\int_{Q\setminus\SC}\!\!\!
\big(\bbC(e(\bar u_{\tau h})){-}\bbC(e(v_h))\big)
{:}e(\bar u_{\tau h}{-}v_h)
\,\d x\d t
\\&\nonumber\le
\int_{Q\setminus\SC}\!\!\!\!\bar F_\tau{\cdot}(\bar u_{\tau h}{-}v_h)
-\big(\bbD e(\DT u_{\tau h})+
\bbC(e(v_h))\big)
{:}e(\bar u_{\tau h}{-}v_h)
\,\d x\d t
\\\displaybreak
\nonumber&+\int_{\SC}\!\!\underline z_{\tau h}\bbA
\JUMP{\bar u_{\tau h}}{\cdot}\JUMP{\bar u_{\tau h}{-}v_h}\,\d S\d t
+\int_{\Snew}\!\bar f_\tau{\cdot}(\bar u_{\tau h}{-}v_h)\,\d S\d t
\\\nonumber&\le\int_{Q\setminus\SC}\!\!\!\bar F_\tau{\cdot}(\bar u_{\tau h}{-}v_h)
-\bbC(e(v_h)){:}e(\bar u_{\tau h}{-}v_h)
-\bbD e(\DT u_{\tau h}){:}e(\bar u_{\tau h}{-}v_h)\,\d x\d t
\\&+\!\int_{\SC}\!\!\!\underline z_{\tau h}\bbA
\JUMP{\bar u_{\tau h}}{\cdot}\JUMP{\bar u_{\tau h}{-}v_h}\,\d S\d t
+\!\int_{\Snew}\!\!\!\bar f_\tau{\cdot}(\bar u_{\tau h}{-}v_h)\,\d S\d t
-\!\int_{Q\setminus\SC}\!\!\!\!\!\bbD e(\dot u_{\tau h}){:}e(\bar u_{\tau h}{-}v_h)\,\d x\d t\ .
\label{strong-conver}\end{align}
The last term can be estimated as follows, which shows its nonpositivity
\begin{align}\label{a-x}
&-\int_{Q\setminus\SC}\!\!\!\!\bbD e(\dot u_{\tau h}){:}e(\bar u_{\tau h}{-}v_h)\,\d x\d t\nonumber
=\sum_{k=1}^{T/\tau}\int_{\Omega\setminus\GC}\!\!\!\!\bbD e(u_{\tau h}^{k-1}{-}u_{\tau h}^k){:}e(u_{\tau h}^k)
+\int_{Q\setminus\SC}\!\!\!\!\bbD e(\dot u_{\tau h}){:}e(v_h)\,\d x\d t
\\\nonumber &\qquad\qquad
\le\sum_{k=1}^{T/\tau}\int_{\Omega\setminus\GC}\!\frac12\bbD e(u_{\tau h}^{k-1}){:}e(u_{\tau h}^{k-1})
-\frac12\bbD e(u_{\tau h}^{k}){:}e(u_{\tau h}^{k})\,\d x
+\int_{Q\setminus\SC}\!\!\!\!\bbD e(\dot u_{\tau h}){:}e(v_h)\,\d x\d t\nonumber 
\\ &\qquad\qquad
=\frac12\int_{\Omega\setminus\SC}\!\!\!\!\bbD e(u_0){:}e(u_0)
-
\bbD e(u_{\tau h}^{T/\tau}){:}e(u_{\tau h}^{T/\tau})\,\d x+
\int_{Q\setminus\SC}\!\!\!\!\bbD e(\dot u_{\tau h}){:}e(v_h)\,\d x\d t  .
\end{align}
Now, using \eqref{ass-mesh}, we take $v_h\to u$ in 
$L^p(I;W^{1,p}(\Omega{\setminus}\GC;\R^d))$ and  notice that the right-hand side in \eqref{a-x} tends to zero which  forces  the left-hand side of the 
expression in \eqref{strong-conver} to vanish, as well. 
This yields strong convergence 
$\bar u_{\tau h}\to u$ in $L^p(I;W^{1,p}(\Omega\setminus\GC;\R^d))$.
By interpolation of $L^\infty(I;W^{1,p}(\Omega\setminus\GC;\R^d))$
and $L^p(I;W^{1,p}(\Omega{\setminus}\GC;\R^d))$, we eventually obtain 
\eqref{conv-u}.


The limit passage in \eqref{disc-u+} is simple with the kinds of convergence proved above. 

For the limit passage in \eqref{disc-z-1+}, we can always consider $\tau>0$ so small, 
that chosen time instants $t_1<t_2$ are  elements of $\{k\tau\}_k$ and we write
\begin{align}\nonumber
&\int_{t_1}^{t_2}\!\calR_1\big(\bar u_{\tau h};\DT z_{\tau h}\big)\,\d t
=\int_{t_1}^{t_2}\int_{\Gamma_C}\!\!
\alpha(\JUMP{\pwc{u}{\tau h}})|\pwl{\DT{z}}{\tau h}|\,\d S\d t
\\&\hspace{4em}
=\int_{t_1}^{t_2}\int_{\Gamma_C}\!\!\alpha\big(\JUMP{u_{\tau h}}\big)|\pwl{\DT{z}}{\tau h}|\,\d
S\d t
+\int_{t_1}^{t_2}\int_{\Gamma_C}\!\!\big(\alpha\big(\JUMP{\pwc{u}{\tau h}}\big)
-\alpha\big(\JUMP{u_{\tau h}}\big)\big)
|\pwl{\DT{z}}{\tau h}|\,\d S\d t. 
\label{7.3}\end{align}
We can see that we need the 
convergence of $\alpha(\JUMP{u_{\tau h}})\to \alpha(\JUMP{u})$ 
in $C(\overlineSC)$. We observe that the sequence 
$\{u_{\tau h}:I\to H^1(\Omega{\setminus}\GC)\}_{\tau>0,h>0}$ 
is equicontinuous because
\begin{align}\nonumber
\big\|u_{\tau h}(t_1){-}u_{\tau h}(t_2)\big\|_{H^1(\Omega{\setminus}\GC)}\le
\Big\|\int_{t_1}^{t_2}\!\!\!\DT u_{\tau h}\,\d t\Big\|_{H^1(\Omega{\setminus}\GC)}
\le\int_{t_1}^{t_2}\!\big\|\DT u_{\tau h}\big\|_{H^1(\Omega{\setminus}\GC)}\d t
\\\label{AC}
\le\big\|\DT u_{\tau h}\big\|_{L^2(I;H^1(\Omega{\setminus}\GC))}
\|1\|_{L^2([t_1,t_2])}
=|t_1{-}t_2|^{1/2}\big\|\DT u_{\tau h}\big\|_{L^2(I;H^1(\Omega{\setminus}\GC))}
\end{align}
for any $0\le t_1<t_2\le T$. By compactness of the embedding 
$W^{1,p}(\Omega{\setminus}\GC)\subset
W^{1-\epsilon,p}(\Omega{\setminus}\GC)$ and 
by Arzel\'a-Ascoli-type arguments based on \eqref{AC}, 
cf.\ \cite[Lemma~7.10]{NPDE_roubicek},
we have the strong convergence $u_{\tau h}\to u$ in 
$C(\bar I;W^{1-\epsilon,p}(\Omega{\setminus}\GC))$ with a small $\epsilon>0$. 
By the trace operator $W^{1-\epsilon,p}(\Omega{\setminus}\GC)\to C(\overlineGC)$,
we then have also the convergence $\JUMP{u_{\tau h}}\to\JUMP{u}$
in $C(\overlineSC)$. Using \eqref{anot1}, we eventually have
$\alpha(\JUMP{u_{\tau h}})\to\alpha(\JUMP{u})$ in $C(\overlineSC)$. 
This allows us to pass to the limit in the first term of \eqref{7.3}. 
Since $|\pwl{\DT{z}}{\tau h}|\to |\pwl{\DT{z}}{}|$ weakly* in the sense 
of measures on $\overlineSC$, we conclude that
\begin{equation}
\label{measure-convergence}
\alpha\big(\JUMP{u_{\tau h}}\big)|\pwl{\DT{z}}{\tau h}|\to
\alpha\big(\JUMP{u}\big)|\pwl{\DT{z}}{}|\qquad \text{ weakly* in the sense
of measures on $\overlineSC$,}
\end{equation}
where $|\DT z|$ denotes the
variation of the measure $\DT z$; in fact, here simply $|\DT z|=-\DT
z$, since $\DT z\le0$. Using again \eqref{AC}, we have 
\begin{align}
\big\|\bar u_{\tau h}{-}u_{\tau h}\big\|_{L^\infty(I;H^1(\Omega{\setminus}\GC;\R^d))}
\le\tau^{1/2}\big\|\DT u_{\tau h}\big\|_{L^2(I;H^1(\Omega{\setminus}\GC;\R^d))}
\end{align}
and, by interpolation with the 
$L^\infty(I;W^{1,p}(\Omega{\setminus}\GC;\R^d))$-estimate, we also obtain that for 
$\epsilon>0$ arbitrarily small $\big\|\bar u_{\tau h}{-}u_{\tau h}
\big\|_{L^\infty(I;W^{1,p-\epsilon}(\Omega{\setminus}\GC;\R^d))}\to0$. For $\epsilon<p{-}d$, 
we have still the continuous trace operator $W^{1,p-\epsilon}(\Omega{\setminus}\GC;\R^d)
\to C(\overlineGC)$. Using \eqref{anot1}, we also have
$\alpha(\JUMP{\bar u_{\tau h}})-\alpha(\JUMP{u_{\tau h}})\to0$ in 
$L^\infty(I;C(\overlineGC))$. Since $\{\pwl{\DT{z}}{\tau h}\}_{\tau>0,h>0}$ is 
bounded in $L^1(\SC)$, we then conclude that the  second term on the right-hand
side of \eqref{7.3} tends to zero as $\tau \to 0$. Thus, altogether 
$\int_{t_1}^{t_2}\!\calR_1(\bar u_{\tau h};\DT z_{\tau h})\,\d t
\to\int_{t_1}^{t_2}\!\calR_1(u;\DT z)\,\d t$.
Furthermore, by weak lower semicontinuity, we have
$\liminf_{\tau\to0,\,h\to0}\int_{t_1}^{t_2}\!\calR_2(\DT u_{\tau h})\,\d t
\ge\int_{t_1}^{t_2}\!\calR_2(\DT u)\,\d t$.
If $t_1\in [0,T)$ is such that  $u_{\tau h}(t_1)\to u(t_1)$ strongly in 
$W^{1,p}(\Omega\setminus\Gamma_C;\R^d)$ then the remaining terms in 
\eqref{disc-z-1+} can be handled even by semicontinuity and continuity 
to obtain (notice that $\varphi$ is $p$-Lipschitz continuous) 
\begin{align}
\Phi\big(u(t_2),z(t_2)\big)+\int_{t_1}^{t_2}\!\!\calD(u;\DT u,\DT z)\,\d t
\le
\Phi\big(u(t_1),z(t_1)\big)
+\int_{t_1}^{t_2}\!\int_\Omega\!\FRM{\cdot}\DT u\,\d x\d t
+\int_{t_1}^{t_2}\!\int_{\GN}\!\!\fRM{\cdot}\DT u\,\d S \d t\,.
 \label{total-energy0}\end{align}

The limit passage in \eqref{disc-z-2+} needs an explicit construction
of a so-called mutual recovery sequence. Take $t\in[0,T]$ for which 
\eqref{disc-z-2+} holds. We want to prove \eqref{semistab}. Take arbitrary 
$\tilde z\in L^\infty(\GC)$ such that $\tilde z\le z(t,\cdot)$ because 
otherwise the inequality holds trivially. Here we can take 
\begin{align}
\tilde z_{\tau h}(t,x)=\begin{cases}
\bar z_{\tau h}(t,x)\Pi_h\Big(\frac{\tilde z(x)}{z(t,x)}\Big)
&\mbox{ if $z(t,x)\ne 0$}\\0 & \mbox{ otherwise,}
\end{cases} 
\end{align}
where $\Pi_h$ is the projector $L^\infty(\GC)\to Z_h$ making piecewise
constant averages, cf.\ \cite[Sect.4.6]{MieRou09NARI}. We see that 
$\tilde z_{\tau h}\le z_{\tau h}$. Notice that $\Pi_h(\tilde z/ z)\to \tilde z/z$ 
as $h\to 0$ in any $L^p(\GC)$ with $+\infty >p\ge 1$. Hence, 
$\tilde z_{\tau h}\to \tilde z$ weakly* in $L^\infty(\GC)$. We get from \eqref{disc-z-2+} that
\begin{align}
0&\le \lim_{\tau\to 0}\lim_{h\to 0}(\Phi(\bar u_{\tau h}(t),\tilde z_{\tau h}) 
+ \mathcal{R}_1(\bar u_{\tau h}(t),\tilde z_{\tau h}-\bar z_{\tau h}(t))
-\Phi(\bar u_{\tau h}(t),\bar z_{\tau h}))\nonumber\\
&= \lim_{\tau\to 0}\lim_{h\to 0}\int_{\GC}\left(\frac12\dela\JUMP{\bar u_{\tau h}}{\cdot}\JUMP{\bar u_{\tau h}}
+\alpha\JUMP{\bar u_{\tau h}}\right)(\tilde z_{\tau h}-\bar z_{\tau h}(t))\,\d S\nonumber\\
&=\int_{\GC}\!\Big(\frac12\dela\JUMP{u}{\cdot}\JUMP{u}+\alpha\JUMP{u}\Big)(\tilde z{-}z(t))\,\d S\nonumber
=\Phi(u(t),\tilde z)+\mathcal{R}_1(u(t),\tilde z{-}z(t))-\Phi(u(t),z(t))\,,
\end{align}  
which proves \eqref{semistab}.

As the considered finite-element discretization always exists
as well as solutions to the recursive problem \eqref{semi-impl+},
the existence of an energetic solution is proved by this
constructive method.
$\hfill\Box$

\bigskip

\begin{remark}[{\it Energy conservation}]
\upshape
The equality in \eqref{total-energy} is still unclear. Let us outline the
difficulties. Within this model, the proof would need \eqref{posit} to hold for
$p=2$, which however contradicts our assumption $p>d$ in the physically
interesting multidimensional cases.
Assume further also that \eqref{ass-monot-phi} hold with $p=2$, so that  we also have  
\begin{align}\label{stress-growth}
&\exists C>0\ \forall e\in\R_{\mathrm{sym}}^{d\times d}:\qquad |\varphi'(e)|\le C(1+|e|)\ .\end{align}
Then we can prove the energy inequality opposite to \eqref{total-energy}. Consider for simplicity 
$t_1=0$ and  $t_2=T$. Following  \cite{rr+tr2} we notice that 
\begin{align}\label{subdifferential}
\Phi(u(T),z(T))-\Phi(u_0,z_0) +\int_{\overlineSC}\!\!
\alpha(\JUMP{u})|\DT z|\,\d S\d t
\ge\int_0^T\!\!\big\langle\lambda,\DT u\big\rangle\,\d t\ ,
\end{align}  
where $\lambda\in L^2(0,T;W^{1,p}(\Omega{\setminus}\GC;\R^d)^*)$ 
with $\lambda(t)\in\partial_u\Phi(u(t),z(t))$ for a.a.~$t\in(0,T)$
 and where ``$|\DT z|\,\d S\d t$'' denotes the integration with 
respect to the Borel measure $|\DT z|$ on $\overlineSC$
.
We have that $\lambda\in\partial_u\Phi(u,z)$ if and only if there 
is $\ell\in\partial\mathcal{J}(u)$ and for all $v\in W^{1,p}(\Omega\setminus\GC;\R^d)$
$$
\big\langle\lambda,v\big\rangle=\int_\Omega\mathbb{C}(e(u)){:}e(v)\,\d x+\int_{\GC}z\mathbb{A}\JUMP{u}{\cdot}\JUMP{v}\,\d S +\big\langle\ell,v\big\rangle\ ,$$
where $\langle\cdot,\cdot\rangle$ denotes the appropriate duality pairing.
Here the convex function $\mathcal{J}: H^1(\Omega\setminus\GC;\R^d)\to [0;+\infty]$ is defined as follows
$$
\mathcal{J}(u)=\begin{cases}
0 &\mbox {if } \JUMP{u}{\cdot}
\norm\ge 0\ ,\\
+\infty &\mbox{ otherwise.} 
\end{cases}
$$
In order to prove \eqref{subdifferential} for a chosen selection 
$\lambda(t)\in\partial_u\Phi(u(t),z(t))$  we use \eqref{semistab} and approximations of 
 Stieltjes integrals by Riemann sums; \cite{DMFraToa05}. Take $n\in\N$  and a 
partition $0=t^n_0<t^n_1<\ldots <t^n_{N_n}=T$ such that $\max_i(t^n_{i+1}-t^n_{i})<1/n$ 
and such that the function $\mathcal{A}_n:[0,T]\to L^\infty(\GC)$ defined as 
$\mathcal{A}_n(t)=\alpha(\JUMP{u(t^n_{i-1})})$ if $t^n_{i-1}<t\le t^n_{i}$. We can assume, without loss 
of generality, that \eqref{semistab} holds for all timesteps $\{t^n_i\}_{i=0}^{N_n-1}$.  
 The sequence $\{\mathcal{A}_n\}$ satisfies
 \begin{align}\label{passage}
\mathcal{A}_n\to \alpha(\JUMP{u})\mbox{ in }L^\infty(\GC) \mbox{ as } n\to\infty\ .
\end{align}
  This limit passage follows from the fact that $u\in H^1([0,T]; H^1(\Omega\setminus\GC;\R^d))$, 
so that $u:[0,T]\to H^1(\Omega\setminus\GC;\R^d)$ is H\"{o}lder continuous on $[0,T]$, 
which implies uniform continuity of $\JUMP{u}:[0,T]\to L^\infty(\GC)$. Finally, \eqref{passage} 
follows by approximations of uniformly continuous functions by piecewise constant interpolants. 
Further, we test \eqref{semistab} at $t^n_{i-1}$ against $\tilde z=z(t^n_i)$ to get
 \begin{align}\label{semistab-ex}
 &\Phi(u(t^n_{i-1}),z(t^n_{i-1}))\le \Phi((u(t^n_{i-1}),z(t^n_i)) 
+\mathcal{R}_1(u(t^n_{i-1}),z(t^n_i)-z(t^n_{i-1}))\nonumber\\&
 =\Phi((u(t^n_{i}),z(t^n_i))-\int_{\GC}\!\mathcal{A}_n(t)(z(t^n_i)-z(t^n_{i-1}))\,\d S
-\int_{t^n_{i-1}}^{t^n_i}\!\!\big\langle\lambda_n(s), \DT u(s)\big\rangle\, \d s
 \end{align}  
whenever $\lambda_n(s)\in\partial_u\Phi(u(s),z(t^n_i))$ for any $s\in(t^n_{i-1},t^n_i]$. 
The last equality is based on the chain rule for $\frac{\d}{\d t}\Phi(u(t), z(t^n_{i}))$; 
cf.\cite[Prop.IX.4.11]{visintin96}. Notice that $\langle\lambda_n,\DT u\rangle$ is well-defined 
even for $\DT u \in H^1(\Omega{\setminus}\GC;\R^d)$ due to \eqref{stress-growth}. 
Summing up \eqref{semistab-ex} for $i=1,\ldots, N_n$ and passing to the limit inferior 
for $n\to \infty$ (in fact, the for the dissipative term on the left hand side even the 
limit exists), we have in view of \eqref{passage} and \cite{rr+tr2} that
$$
\Phi(u(T),z(T))-\Phi(u_0,z_0) -\int_{\overlineSC}
\!\!\alpha(\JUMP{u})\DT z\,\d S\d t
\ge\int_0^T\!\big\langle\lambda(s),\DT u(s)\big\rangle\,\d s\ .
$$
Further, notice that $\int_{t^n_1}^{t^n_2}\langle \ell,\DT u\rangle\,\d s
=\mathcal{J}(u(t^n_2))- \mathcal{J}(u(t^n_1))=0$, hence we have
\begin{align}\label{opposite0}
 &\Phi(u(T),z(T))-\Phi(u_0,z_0)-\!\int_{\overlineSC}\!\!\!
\alpha(\JUMP{u})\DT z\,\d S\d t
\ge
\!\int_Q\!
\mathbb{C}(e(u)){:}e(\DT u)\,\d x\d t
+
\!\int_{\SC}\!\!\!
z\mathbb{A}\JUMP{u}{\cdot}\JUMP{\DT u}\,\d S
\,\d t\,.
\end{align}
Taking into account \eqref{e:weak-momentum-var} we see that $\mathfrak{J}$ defined by
\begin{align}
\big\langle\mathfrak{J},v\big\rangle:=\int_{\Omega\setminus\GC}\!\!\!\!
\big(\bbD e(\DT{u}){+}\bbC(e(u))\big){:}e(v)
\,\d x
+\int_{\GC}\!\!\!
z\bbA\JUMP{u}
{\cdot}\JUMP{v}\d S
-\int_\Omega\!\FRM{\cdot} v\,\d x-\int_{\GN}\!\!\fRM{\cdot} v\,\d S
\label{e:weak-momentum-var1}
\end{align}
is such that 
$$
\int_0^T\mathcal{J}(v)\,\d t-\int_0^T\mathcal{J}(u)\,\d t
\ge \int_0^T\big\langle\mathfrak{J},v{-}u\big\rangle\,\d t\ , 
$$
i.e. $\mathfrak{J}(t)\in\partial\mathcal J(u(t))$. Therefore, 
$\int_0^T\langle\mathfrak{J},\DT u\rangle\,\ \d t=0$ and together 
with \eqref{opposite0} we get that 
\begin{eqnarray*}
\Phi\big(u(t),z(t)\big)+\int_0^t\!\calD(u;\DT u,\DT z)\,\d t
\ge
\Phi\big(u_0,z_0\big)
+\int_0^t\!\int_\Omega\!\FRM{\cdot}\DT u\,\d x\d t
+\int_0^t\!\int_{\GN}\!\!\fRM{\cdot}\DT u\,\d S \d t\ .
\end{eqnarray*}
\end{remark}

\begin{remark}[{\it Linear but nonlocal elasticity}]\label{rem-linear}
\upshape 
We can also consider linear but nonlocal stress
\begin{align}
\sigma=\bbD e(\DT u)+\bbC e(u)+
\int_\Omega\bbH(x,\xi)\big(e(u(x){-}u(\xi))\big)\,\dd\xi,
\end{align}
with some (presumably only small) 4th-order tensor $\bbH$. This yields a 
contribution 
$$
\frac12\!\sum_{i,j,k,l=1}^d\int_{\Omega\times\Omega}
\bbH_{ijkl}(x,\xi)e_{ij}(u(x){-}u(\xi))e_{kl}(u(x){-}u(\xi))\,\dd\xi\dd x
$$ 
to the stored energy. Assuming that, for some $\epsilon>0$ and $\zeta>0$ 
and for all $x,\xi\!\in\!\Omega$, it holds 
$\bbH_{ijkl}(x,\xi)\ge\zeta\delta_{ik}\delta_{jl}/|x{-}\xi|^{d+2\epsilon}$, 
we obtain estimates in $W^{1+\epsilon,2}(\Omega{\setminus}\GC;\R^d)$
instead of $W^{1,p}(\Omega{\setminus}\GC;\R^d)$
we used above. For $d=2$, any $0<\epsilon<1$ is sufficient for 
the convergence and, for $\epsilon<1/2$, we can still use P1-elements
for the space discretization. For $d=3$, $\epsilon>1/2$ is needed
for the compact embedding $W^{1+\epsilon,2}(\Omega{\setminus}\GC)
\subset C(\barOmegaone{\cup}\barOmegatwo)$ and thus for the convergence, which
unfortunately excludes discontinuities of gradients and thus P1-element
so that P2-elements would have to be used. In the latter case,
one could also consider the concept of nonsimple materials. A definite benefit would be an energy conservation \cite{rr+tr2}.
\end{remark}

\begin{remark}[{\it Towards the rate-independent evolution}]\label{rem-viscosity}
\upshape 
We can replace $\mathbb{D}$ with $\varepsilon\mathbb{D}$ in the model and, like in
the mixity-insensitive case \cite{Roub13ACVE,RoPaMA13QACV}, study the 
vanishing-viscosity limit 
for $\varepsilon\!\searrow\!0$ towards the rate-independent limit. 
Let us denote $(u_\varepsilon,z_\varepsilon)$ a weak solution to the problem (\ref{adhes-classic}), 
i.e., it satisfies Definition~\ref{def4} with our modified viscosity. Then we get the 
following a-priori estimate 
$\sqrt{\varepsilon}\|e(\DT u_\varepsilon)\|_{L^2(Q\setminus\SC)}\le S_0$. We consider a weak$*$ limit
$(u,z)\in L^\infty(0,T;W^{1,p}(\Omega{\setminus}\GC;\R^d))\times
(L^\infty(\SC)\cap {\rm BV}([0,T];L^1(\GC)))$  of a (sub)sequence of solutions 
for $\varepsilon\!\searrow\!0$. Further, using \eqref{ass-monot-phi}, we estimate for 
$\varepsilon\!\searrow\!0$ that:
\begin{align*}
\eta\big\|e(u_\varepsilon{-}u)\big\|_{L^p(Q;\R^{d\times d})}^p&\le
\int_{Q\setminus\SC}\!\!\!\big(\mathbb{C}(e(u_\varepsilon))-\mathbb{C}(e(u))\big)
{:}e(u_\varepsilon{-}u)\,\d x\d t\\
&\le \int_{Q\setminus\SC}\!\!\!(\mathbb{C}(e(u_\varepsilon))-\mathbb{C}(e(u)))
{:}e(u_\varepsilon{-}u)\,\d x\d t\\[-.5em]
&\qquad\quad+\int_{\SC}\!\!\!z_\varepsilon\mathbb{A}\JUMP{u_\varepsilon{-}u}
{\cdot}\JUMP{u_\varepsilon{-}u}\,\d S\d t\\
&\le-\int_{Q\setminus\SC}\!\!\!\varepsilon\mathbb{D}e(\DT u_\varepsilon){:}e(u_\varepsilon{-}u)
+\mathbb{C}(e(u)){:}e(u_\varepsilon{-}u)\,\d x\d t\\[-.5em]
&\qquad\qquad-\int_{\SC}\!\!\!z_\varepsilon\mathbb{A}\JUMP{u}{\cdot}\JUMP{u_\varepsilon{-}u}
+(z_\varepsilon{-}z)\mathbb{A}\JUMP{u_\varepsilon}{\cdot}\JUMP{u_\varepsilon}\,\d x\d t\\[-.5em]
&\qquad\qquad\qquad
+\int_Q \!F{\cdot}(u_\varepsilon{-}u)\,\d x\d t
+\int_{\Snew}\!\!f{\cdot}(u_\varepsilon{-}u)\,\d S\d t\to0\ .
\end{align*}
Hence, we have for almost every $t\in[0,T]$ that $e(u_\varepsilon(t))\to e(u(t))$ in 
$L^p(\Omega{\setminus}\GC;\R^{d\times d})$. We thus obtain an approximable solution
in the sense of \cite{DDMM08VVAQ,KnMiZa08ILMC,ToaZan09AVAQ} to the rate-independent 
mixity-sensitive model. Moreover, also  
(\ref{semistab}) holds for a.a. $t\in[0,T]$. It is not obvious, however, if it holds that 
$\liminf_{\varepsilon\to 0}\int_{t_1}^{t_2}\!\int_{\overlineGC}\alpha(\JUMP{u_\varepsilon})\DT z_\varepsilon\,\d S\d t
\ge\int_{t_1}^{t_2}\!\int_{\overlineGC}\alpha(\JUMP{u})\DT z\,\d S\d t$, which is needed to pass to the 
limit in the energy (im)balance. In this sense, the associative model
from \cite{tr+mk+jz,RoMaPa13QMMD,RoPaMA13QACV,RoPaMa??LSAQ}
remains still more justified for the merely rate-independent mixity-sensitive evolution.

\end{remark}
\subsection*{Acknowledgments}
The authors are deeply indebted to prof.\ Vladislav Manti\v c for valuable comments.   
This contribution has been partly supported by 
the Junta de Andaluc\'{\i}a (Proyecto de Excelencia 
P08-TEP-4051) and the Spanish Ministry of Education and Science (Proyecto
MAT2009 - 140022), as well as by the grants 
201/09/0917, 201/10/0357, 201/12/0671, and 13-18652S (GA \v CR),
and the institutional support RVO:61388998 (\v CR).




{\small

\begin{thebibliography}{SK}

\normalsize \baselineskip=12pt


\vspace{-.5em}\bibitem{BanAsk00NFCI} {\sc L. Banks-Sills, D. Askenazi}:
{A note on fracture criteria for interface fracture}. {\it Intl. J.
Fracture}  {\bf 103} (2000), 177-188.




\vspace{-.5em}\bibitem{BBR1} {\sc E.Bonetti, G.Bonfanti,  R.Rossi}:
\newblock  Global existence for a contact problem with adhesion.
\newblock  {\it Math. Methods Appl. Sci.} {\bf 31} (2008),
1029--1064.
%
%
%
%

\vspace{-.5em}\bibitem{CaDaMo03NSMM}
{\sc P. P. Camanho, C.G. D\'avila, M.F. De Moura}:
Numerical simulation of mixed-mode progressive
delamination in composite materials.
{\it J. Composite Mater.} {\bf 37} (2003), 1415-1438.

\vspace{-.5em}\bibitem{CoCa11ES}
{\sc P. Cornetti, A.  Carpinteri}:
Modelling the {FRP}-concrete delamination
by means of an exponential softening law,
{\it Engineering Structures} {\bf 33} (2011), 1988-2001.

\vspace{-.5em}\bibitem{DMFraToa05} {\sc G. Dal Maso, G. Francfort,  R.
Toader}: \newblock Quasistatic crack growth in nonlinear elasticity.
\newblock  {\it Arch. Rational Mech. Anal.} {\bf 176} (2005),
165--225.


\vspace{-.5em}\bibitem{ERDC90FEBI} {\sc A.G. Evans, M. R\"{u}hle,
 B.J. Dalgleish, P.G. Charalambides}: {The fracture energy of bimaterial
interfaces}. {\it Metallurgical Transactions A} {\bf 21A} (1990),
2419-2429.



\vspace{-.5em}\bibitem{DDMM08VVAQ}
{\sc G. Dal Maso, A. DeSimone, M.G. Mora, M. Morini}:
A vanishing viscosity approach to quasistatic evolution 
                  in plasticity with softening, 
{\it Arch. Rational Mech. Anal.} {\bf 189}
(2008), 469-544.


\vspace{-.5em}\bibitem{Fre82} {\sc M. Fr\'emond}: Equilibre des
structures qui adh\`erent \`a leur support. 
{\it Comptes Rendus
Acad. Sci. Paris} \textbf{295} (1982), 913--916.

\vspace{-.5em}\bibitem{Fre87} {\sc M. Fr\'emond}: Adh\`erence
des solides. 
{\it J. M\'ech. Th\'eorique Appliq.} 
\textbf{6} (1987), 383--407.



\vspace{-.5em}\bibitem{HutSuo92MMCL} {\sc J.W. Hutchinson, T.Y. Wu}:
Mixed mode cracking in layered materials. In:
{\it Advances in Appl.\ Mech.}
(Eds. J.W. Hutchinson, T.Y.Wu), Acad. Press, New
York, 1992, pp.63-191.

\vspace{-.5em}\bibitem{KoMiRo} {\sc M. Ko\v cvara, A. Mielke,  T.
Roub\'\i\v cek}:  \newblock  A rate-independent approach to the
delamination problem. \newblock
{\it Mathematics and Mechanics of Solids} {\bf 11} (2006), 423--447.
 
\vspace{-.5em}\bibitem{KnMiZa08ILMC}
{\sc D. Knees, A. Mielke, C. Zanini}:
On the inviscid limit of a model for crack propagation.
{\it Math. Models Meth. Appl. Sci.} {\bf 18} (2008),
1529-1569.

\vspace{-.5em}\bibitem{LaOrSu10ESRM} 
{\sc C.J. Larsen, C. Ortner, E. S\"uli}: 
Existence of solution to a regularized model of dynamic
fracture. {\it Math. Models Meth. Appl. Sci.} {\bf 20} (2010),
1021-1048.

\vspace{-.5em}\bibitem{LieCha92ASIF} {\sc K.M. Liechti, Y.S. Chai}:
{Asymmetric shielding in interfacial fracture under in-plane shear},
{\it J. Appl. Mech.} {\bf 59} (1992), 295-304.

\vspace{-.5em}\bibitem{Mant08DRLM} {\sc V. Manti\v{c}}: Discussion on
the reference lenght and mode mixity for a bimaterial interface.
{\it J. Engr. Mater. Technology} {\bf 130} (2008), 045501-1-2.

\vspace{-.5em}\bibitem{Miel05ERIS} {\sc A. Mielke}: 
Evolution in rate-independent systems.
In: {\it Handbook of Differential Equations,
Evolutionary Equations, 2} (Eds.: Dafermos, C.M., Feireisl, E.),
Elsevier, Amsterdam, 2005, pp.~461--559.

\vspace{-.5em}\bibitem{Miel11altDEMF}
{\sc A. Mielke}: Differential, energetic and metric formulations for 
rate-independent processes. In: {\it Nonlinear PDE's and Applications} 
(Eds.: L.Ambrosio, G.Savar\'e), Springer, 2011, pp.~87--170.

\vspace{-.5em}\bibitem{MieRou09NARI}
 {\sc A. Mielke, T. Roub{\'\i}{\v{c}}ek}:
Numerical approaches to rate-independent processes and applications in 
inelasticity. {\it Math. Model. Numer. Anal.} {\bf 43} (2009), 399--428.



\vspace{-.5em}\bibitem{MieThe04RIHM} {\sc A. Mielke, F. Theil}: 
On rate-independent hysteresis models.
{\it Nonlin. Diff. Equations Appl.} \textbf{11} (2004),
151--189. 

\vspace{-.5em}\bibitem{MiThLe02VFRI} 
{\sc A. Mielke, F. Theil, V. Levitas}: A variational formulation of 
rate-independent phase transformations using an extremum principle. 
\newblock
{\it Arch.\ Rational Mech.\ Anal.} \textbf{162} (2002), 137--177.

\vspace{-.5em}\bibitem{Mukh01BIEC}
{\sc S. Mukherjee}: On boundary integral equations for cracked and for 
thin bodies. {\it Mathematics and Mechanics of Solids}, {\bf 6}
(2001), 47-64.

\vspace{-.5em}\bibitem{ParCan97BEMF}
{\sc F. Par\'{\i}s, J. Ca\~nas}:
{\it Boundary Element Method, Fundamentals and Applications}
Oxford Univ.\ Press, Oxford, 1997.







\vspace{-.5em}\bibitem{rr+tr} {\sc R. Rossi, T. Roub\'\i\v cek}:
Thermodynamics and analysis of rate-independent adhesive contact at
small strains.
{\it Nonlin. Anal.} {\bf 74} (2011), 3159-3190.

\vspace{-.5em}\bibitem{rr+tr2} {\sc R. Rossi, T. Roub\'\i\v cek}:
Adhesive contact delaminating at mixed mode, its thermodynamics and analysis. 
(Preprints arXiv:1110.2794). {\it Interfaces and Free Boundaries} 
{\bf 14} (2013), in print. DOI: 10.4171/IFB/293

\vspace{-.5em}\bibitem{NPDE_roubicek} {\sc T. Roub\'\i\v cek}: {\it
Nonlinear Partial Differential Equations with Applications}.
Birkh\"auser, Basel, 
(2nd.\,ed.), 2013.


\vspace{-.5em}\bibitem{tr2} {\sc T. Roub\'\i\v cek}: Rate independent
processes in viscous solids at small strains. {\it Math. Methods
Appl. Sci.} {\bf 32} (2009), 825--862.




\vspace{-.5em}\bibitem{Roub13ACVE}
{\sc T. Roub\'\i\v cek}: Adhesive contact of visco-elastic bodies
and defect measures arising by vanishing viscosity.
{\it SIAM J. Math. Anal.} {\bf 45} (2013), 101-126.

\vspace{-.5em}\bibitem{Roub??MDLS}
{\sc T. Roub\'\i\v cek}: Maximally-dissipative local solutions to rate-independent systems
and application to damage and delamination problems. {\it Nonlinear Anal., Th. Meth. Appl.},
submitted. 

\vspace{-.5em}\bibitem{tr+mk+jz} {\sc T. Roub\'\i\v cek, M.Kru\v z\'\i k,
J.Zeman}: Delamination and adhesive contact models and their
mathematical analysis and numerical treatment. Chap.\,9 in: {\it
Math.\ Methods and Models in Composites} (Ed.V.Manti\v c), Imperial
College Press, 2013, pp.349-400. In print, ISBN 978-1-84816-784-1.

\vspace{-.5em}\bibitem{RoMaPa13QMMD} {\sc T. Roub\'\i\v cek, V. Manti\v c,
C.G. Panagiotopoulos}: Quasistatic mixed-mode delamination model.
{\it Disc. Cont. Dynam. Syst., Ser.S} {\bf 6} (2013), 591-610.

\vspace{-.5em}\bibitem{RoPaMA13QACV}
{\sc T. Roub\'\i\v cek, C.G. Panagiotopoulos, V. Manti\v c}:
 Quasistatic adhesive contact of visco-elastic bodies and its numerical treatment 
for very small viscosity. 
{\it Zeitschrift angew. Math. Mech.} {\bf 93} (2013), in print. DOI: 10.1002/zamm.201200239

\vspace{-.5em}\bibitem{RoPaMa??LSAQ}
{\sc T. Roub\'\i\v cek, C.G. Panagiotopoulos, V. Manti\v c}:
Local-solution approach to quasistatic rate-independent mixed-mode delamination.
{\it Math. Models Methods in the Appl. Sci.}, submitted.

\vspace{-.5em}\bibitem{tr-LS-CZ} {\sc T. Roub\'\i\v cek, L. Scardia,  C.
Zanini}: Quasistatic delamination problem.
{\it Cont. Mech. Thermodynam.} {\bf 21} (2009), 223-235.

\vspace{-.5em}\bibitem{tr-os-rv}
{\sc T. Roub\'\i\v cek, O. Sou\v cek, R. Vodi\v cka}:
A model of rupturing lithospheric faults with re-occurring earthquakes.
(Preprint No.2013-013, Ne\v{c}as center, Prague) {\it SIAM J. Appl. Math.}, 
accepted.



\vspace{-.5em}\bibitem{SwLiLo99ITAM} {\sc J.G. Swadener, K.M. Liechti,
A.L. de{L}ozanne}: {The intrinsic toughness and adhesion mechanism
of a glass/epoxy interface}. {\it J. Mech. Phys. Solids} {\bf 47}
(1999), 223--258.

\vspace{-.5em}\bibitem{TMGCP10ACTA} {\sc L. T\'{a}vara, V. Manti\v{c},
E. Graciani, J. Canas,  F. Par\'{\i}s}:
Analysis of a crack in a thin adhesive layer between
orthotropic materials: an application to composite
interlaminar fracture toughness test.
{\it CMES} {\bf 58} (2010), 247-270.

\vspace{-.5em}\bibitem{TMGP11BEMA} {\sc L. T\'{a}vara, V. Manti\v{c}, E.
Graciani, F. Par\'{\i}s}: BEM analysis of crack onset and
propagation along fiber-matrix interface under transverse tension
using a linear elastic-brittle interface model. {\it Engr. Analysis with
Boundary Elements} {\bf 35} (2011), 207--222.

\vspace{-.5em}\bibitem{Thom10PhD}
{\sc M. Thomas}: Rate-independent damage processes in nonlinearly elastic materials.
PhD-thesis, Institut f\"ur Mathematik, Humboldt-Universit\"at zu Berlin, 2010.

\vspace{-.5em}\bibitem{ThoMie10DNEM}
{\sc M. Thomas, A.Mielke}: Damage of nonlinearly elastic materials at small strain 
- existence and regularity results. {\it Z.\ angew.\ Math.\ Mech.} {\bf 90} (2010),
88--112. 

\vspace{-.5em}\bibitem{toupin} {\sc R.A. Toupin}: Elastic materials with
couple stresses. {\it Arch.\ Rational Mech.\ Anal.} {\bf 11} (1962),
385--414.

\vspace{-.5em}\bibitem{ToaZan09AVAQ}
{\sc R.~Toader, C.~Zanini}:
An artificial viscosity approach to quasistatic crack growth.
{\it Boll. Unione Matem. Ital.} {\bf 2} (2009), 1--36.

\vspace{-.5em}\bibitem{TveHut93IPMM} {\sc W. Tveergard, J.V. Hutchinson}: 
{The influence of plasticity on mixed mode interface
toughness}. {\it Intl. J. Physics of Solids} {\bf 41} (1993),
1119-1135.

\vspace{-.5em}\bibitem{visintin96} 
{\sc A. Visintin}: {\it Models of Phase Transitions}. Birkh\"auser, 
Boston, 1996.

\end{thebibliography}

}
\end{document}